\documentclass[preprint,review,colorinlistoftodos]{elsarticle}
\usepackage[margin=1in]{geometry}
\usepackage[hidelinks]{hyperref}
\usepackage{graphicx}
\usepackage{amsmath,amsthm, amssymb, float, bm}
\usepackage[labelfont=bf,size={footnotesize}]{caption}
\usepackage{subcaption}
\usepackage{algorithm}
\usepackage{algpseudocode}
\biboptions{sort&compress}

\usepackage{setspace}
\ifdefined\usetodonotes
\usepackage{todonotes}

\fi
\usepackage{soul}
\usepackage{tcolorbox}
\definecolor{C0}{HTML}{1F77B4}
\definecolor{C1}{HTML}{FF7F0E}
\definecolor{C2}{HTML}{2ca02c}
\definecolor{C3}{HTML}{d62728}
\definecolor{C4}{HTML}{9467bd}
\definecolor{C5}{HTML}{8c564b}

\usepackage{xpatch,environ}
\RenewEnviron{abstract}{%
  \xdef\theabstracttext{%
    \unexpanded{%
      \def\baselinestretch{1}\noindent\unskip\textbf{Abstract}\par\medskip
      \noindent\unskip\ignorespaces}%
    \unexpanded\expandafter{\BODY}%
  }%
}
\def\theabstracttext{}
\xpatchcmd{\pprintMaketitle}
{\ifvoid\absbox}
{\ifx\theabstracttext\empty\else\printtheabstracttext\fi\ifvoid\absbox}
{}{}
\newcommand{\printtheabstracttext}{{%
    \begin{trivlist}
      \normalfont\normalsize
    \item\relax
      \doublespacing\theabstracttext
    \end{trivlist}
  }}

\makeatletter
\def\ps@pprintTitle{%
  \let\@oddhead\@empty
  \let\@evenhead\@empty
  \def\@oddfoot{\reset@font\hfil\thepage\hfil}
  \let\@evenfoot\@oddfoot
}
\makeatother

\usepackage[noabbrev]{cleveref}

\newtheorem{theorem}{Theorem}

\title{Robust, strong form mechanics on an adaptive structured grid: efficiently solving variable-geometry near-singular problems with diffuse interfaces}

\author[au]{Vinamra Agrawal}
\author[uccs]{\texorpdfstring{Brandon Runnels\corref{cor1}}{Brandon Runnels}}
\cortext[cor1]{Corresponding author}
\address[au]{Department of Aerospace Engineering, Auburn University, Auburn, AL USA}
\address[uccs]{Department of Mechanical and Aerospace Engineering, University of Colorado, Colorado Springs, CO USA}

\begin{document}

\begin{abstract}
  Many solid mechanics problems on complex geometries are conventionally solved using discrete boundary methods.
  However, such an approach can be cumbersome for problems involving evolving domain boundaries due to the need to track boundaries and constant remeshing. 
  The purpose of this work is to present a comprehensive strategy for efficiently solving such problems on an adaptive structured grid, while expositing some of the basic yet important nuances associated with solving near-singular problems in strong form.
  We employ a robust smooth boundary method (SBM) that represents complex geometry implicitly, in a larger and simpler computational domain, as the support of a smooth indicator function.
  We present the resulting semidefinite equations for mechanical equilibrium, in which inhomogeneous boundary conditions are replaced by source terms.
  In this work, we present a computational strategy for efficiently solving near-singular SBM-based solid mechanics problems.
  We use the block-structured adaptive mesh refinement (BSAMR) method, coupled with a geometric multigrid solver for an efficient solution of mechanical equilibrium.
  We discuss some of the practical numerical strategies for implementing this method, notably including the importance of grid versus node-centered fields.
  We demonstrate the solver's accuracy and performance for three representative examples: a) plastic strain evolution around a void, b) crack nucleation and propagation in brittle materials, and c) structural topology optimization. 
  In each case, we show that very good convergence of the solver is achieved, even with large near-singular areas, and that any convergence issues arise from other complexities, such as stress concentrations.
\end{abstract}

\ifdefined\usetodonotes
\begin{tcolorbox}[title={\bf Todonote convention}]
  Please use the following convention when making notes:
  \begin{center}\verb-\yourname{Addresseename, we need to XYZ}-\end{center}
  In other words use your name macro to make any comments, then address specific people in the text.
  For instance, if Brandon wants Bob to run more results, he should express this in the following way.
  \begin{center}
    \verb-\brandon{Bob, we need to run more results}-
  \end{center}
\end{tcolorbox}
\listoftodos
\setcounter{page}{0}
\fi
\maketitle

\section{Introduction}
Many computational mechanics problems involve analyzing mechanical systems with highly variable geometry.
Such problems require that the mechanical deformations, and resulting stresses, be resolved subject to a set of complex, time-varying, and sometimes unknown topologies.
Such examples include, but are not limited to, fracture mechanics, problems involving material growth or removal (such as dendrite growth), or structural topology optimization.
In all of these, it is essential to accurately solve mechanical equilibrium equations.
Computational mechanics has historically been overwhelmingly dominated by the finite element method (FEM) due to its ability to conform to arbitrary geometry through iso-parametric elements.
Indeed, FEM is nearly synonymous with computational elasticity.
However, in the case of variable topology, the key advantage of FEM - conformal meshing with isoparametric elements - is less beneficial.
This may necessitate costly mid-simulation remeshing, the use of an explicitly meshed and overlayed boundary, or the use of excessive refinement in anticipation of topological change.

The strong form method with finite differences is an attractive alternative for such problems.
Recently, it was shown by the authors that this method may be coupled to the block-structured adaptive mesh refinement (BSAMR) method, along with the geometric multigrid method, to produce a highly efficient linear elastic solver \cite{runnels2021massively}.
The solver has been applied to numerous small strain \cite{celestine2020experimental,runnels2020phase,gokuli2021multiphase} and finite deformation \cite{strutton2022interface} mechanics problems.
Phase field methods have also been implemented using this method, albeit with problematically slow convergence rates due to limitations of the method that will be addressed in this work \cite{agrawal2021block}.
Many of the problems of interest are reducible to representative volume elements (RVE), for which the finite difference method is ideally suited.
We seek to apply this method to problems in which the geometry may be considered to be variable.

We let the ``smooth boundary method'' (SBM) refer to the approach in which a complex geometry is defined within a simpler computational domain as the support of some smooth indicator function $\phi>0.5$.
Smoothness requires that the transition from solid to void is continuous, and we assume in general that the $\phi$ varies smoothly from $0$ to $1$ over some finite interval.
We consider SBM to refer specifically to the technique of replacing discrete boundary conditions with equivalent source terms, such that the boundary conditions are recovered exactly in the sharp interface limit.
By thus embedding the complex geometry through a diffused interface, SBM circumvents the challenges associated with domain meshing encountered in discrete interface approaches.
The SBM has been used to solve partial differential equations with general boundary conditions on complex boundaries and can easily be coupled with diffuse boundary methods (such as phase field) by evolving the order parameters using a thermodynamic equation.
Some examples of SBM applications include the use of phase field methods to study corrosion in Mg alloys \cite{chadwick2018numerical}, mass flux boundary conditions in fluids \cite{schmidt2022self}, and general partial differential equations \cite{yu2012extended,li2009solving}.

The SBM's efficiency relies heavily on using BSAMR to resolve the diffuse boundary.
When suitably coupled, SBM effectively eliminates the need for explicitly defining the mesh since the interface can be resolved with an appropriate resolution and the mesh can be updated to track the evolving interface.
The above-mentioned BSAMR strategy has been widely used for high-performance computational fluid mechanics problems \cite{zhang2019amrex, hittinger2013block, schornbaum2018extreme, dubey2014survey} and, to some extent, for solid mechanics problems \cite{agrawal2021block, runnels2020phase, runnels2021massively}.
BSAMR stores and evolves each mesh level independently, evolving finer levels with smaller time steps to avoid overly restrictive CFL conditions on coarser levels.
The information between levels is communicated through averaging (fine level to coarse level) and ghost cells (coarse level to fine level). 

Application of SBM in solid mechanics applications can lead to semi-definite problems due to the lack of uniqueness resulting from a mesh-resolved ``void'' region. 
This situation often arises in topology optimization problems where the simulation domain is an output rather than the input but is endemic to any implicit boundary method.
Without properly addressing the semi-definiteness of the operator, the result can be poor (or no) convergence and, worse, an incorrect solution.

In this work, we present computational techniques to allow for the efficient solution of semi-definite smooth boundary problems using the finite difference method with BSAMR.
The paper is structured in the following way.
In \Cref{sec:computational_methods}, the SBM is formalized for elasticity, and the specific challenges of solving semidefinite problems are addressed.
This section also describes some of the computational methods and challenges unique to solving problems of this nature.
In \Cref{sec:examples}, three representative examples are presented that demonstrate the model's effectiveness: (1) plasticity with variable geometry, (2) phase field fracture mechanics, and (3) structural topology optimization.
Each example is somewhat self-contained, so the disparate applications will find relevance in their respective communities.
We conclude by highlighting some limitations of the framework in its current form.

\section{Computational methods}\label{sec:computational_methods}
In this section, we present the key elements of the SBM method and practical strategies for its implementation. 
We first present the formulation of the equations of linear elasticity with traction boundary conditions with SBM. 
We then outline the reflux-free multigrid implementation of the solver using BSAMR.
Next, we emphasize the need for choosing a cell-based indicator field for the stability of the solver. 
Finally, we outline the implementation of material models as a vector space to efficiently work with the solver.

\begin{figure}
  \includegraphics[width=\linewidth]{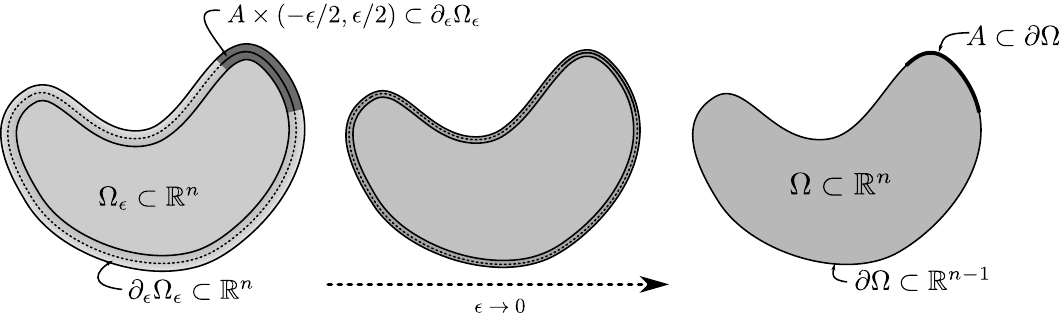}
  \caption{
      Illustration showing the diffuse boundary construction approaching the sharp interface limit.
      The diffuse domain $\Omega_\epsilon$ is the support of $\phi^\epsilon=1$, equal to $\Omega$ in the limit.
      The diffuse boundary $\partial_\epsilon\Omega_\epsilon$ is the $\epsilon$ an open set defined in the region where $0<\phi^\epsilon<1$, and vanishes in the limit.
      It is also defined as the $\epsilon/2$-neighborhood of $\partial\Omega$.
      The shaded region on the left is $A\times(-\epsilon/2,\epsilon/2)$, a subset of $\partial_\epsilon\Omega_\epsilon$ corresponding to an arbitrary subset $A$ of $\partial\Omega$.
  }
\end{figure}

\subsection{Diffuse boundary method for linear elasticity}
In this section we present the diffuse boundary formulation of mechanical equilibrium for a linear elastic material.
We consider a body of interest occupying some region $\Omega\subset\mathbb{R}^3$, with a natural boundary $\partial\Omega$ upon which surface tractions $\bm{t}^0(\bm{y}) \ \ \bm{y}\in\partial\Omega$ are prescribed as boundary conditions.
(We do not consider the diffuse formulation of essential, i.e. displacement, boundary conditions at this time.)
In the discrete setting, the problem is posed as a differential equation with boundary conditions prescribed at the domain boundary.
In the diffuse setting, we represent boundary effects implicitly.
To construct the diffuse problem, we first replace the explicit domain $\Omega$ with a continuous function, called an order parameter $\phi^\epsilon$, that represents $\Omega$.
In the limit as $\epsilon\to0$, the support of $\phi^\epsilon$ is identical to the discrete-boundary domain, $\Omega$; this is called the sharp interface limit.
For $\epsilon>0$, we define the diffuse domain and boundary to be, respectively:
the support of $\phi^\epsilon=1$, denoted $\Omega_\epsilon$;
the support of $\nabla\phi^\epsilon$, denoted $\partial_\epsilon\Omega_\epsilon$.
We require that $\phi^\epsilon=0$ outside of $\partial_\epsilon\Omega_\epsilon\cup\Omega_\epsilon$, and that $\partial_\epsilon\Omega_\epsilon = \partial\Omega\times(-\epsilon/2,\epsilon/2)$, bounding the diffuse region to the $\epsilon$ neighborhood of the discrete boundary.
As long as $\epsilon<<r_{\mathrm{min}}$, the smallest radius of curvature of $\partial\Omega$, we require that there exist some parameterization of $\phi^\epsilon$ within the diffuse boudnary by
\begin{align}
  \phi^\epsilon(\bm{y} + s\bm{n}) = \hat{\phi}(s) \ \ \ \bm{y}\in\partial\Omega, \ \ s\in(-\epsilon/2,\epsilon/2)
\end{align}
where $\bm{n}$ is the normal to the discrete interface.
$\hat{\phi}$ is some Lipschitz function describing the behavior of $\phi$ over the diffuse interface, and is generally a regularized step function over the interval $(-\epsilon/2,\epsilon/2)$.
Its derivative is not defined in the limit as $\epsilon\to0$, but approaches the Dirac delta distribution.
With these definitions and restrictions in hand, it is possible to establish the following theorem, which was presented in \cite{schmidt2022self}:

\begin{theorem}
  Let $\phi^\epsilon$ be an idealized order parameter with length scale $\epsilon$, and let $f$ and $g$ be either scalar or vector-valued bounded functions,
  with $\bm{n}\cdot\nabla g$ bounded in $\partial_\epsilon\Omega_\epsilon$.
  \footnote{We note that the $\bm{n}\cdot\nabla g$ boundeness restriction is erroneously absent from the original presentation of the theorem in \cite{schmidt2022self}.}
  Then the following holds:
  \begin{equation}
    \lim_{\epsilon\rightarrow 0} \int_A \int_{-\epsilon/2}^{\epsilon/2} \big( f \phi^\epsilon + g|\nabla \phi^\epsilon|\big) \, ds\, dA = \int_A g \, dA\quad \forall A \subset \partial\Omega
  \end{equation}
\end{theorem}
We can now present the diffuse interface formulation of mechanical equilibrium.
Recall the usual sharp-interface
equations of momentum conservation in the context of elasticity, with kinematics linearized about some eigenstrain $\bm{\varepsilon}^0$, are:
\begin{subequations}
  \begin{align}
    \mathbb{C}(\bm{x}) \left( \operatorname{grad} \bm{u}(\bm{x}) - \bm{\varepsilon}^0 (\bm{x})\right) - \bm{\sigma}(\bm{x}) &= 0,\;\; \bm{x}\in\Omega & \bm{u}^0(\bm{y}) - \bm{u}(\bm{y}) &= 0,\;\; \bm{y} \in \partial_1\Omega \label{eq:kinematics}\\
    \operatorname{div} \bm{\sigma}(\bm{x}) - \bm{b}(\bm{x}) &= 0, \;\; \bm{x}\in\Omega & \bm{t}^0(\bm{y}) - \bm{\sigma}(\bm{y})\hat{\bm{n}}(\bm{y}) &= 0,\;\;\bm{y}\in\partial_2\Omega.\label{eq:momentum}
  \end{align}
\end{subequations}
\Cref{eq:kinematics,eq:momentum} are the constitutive and mechanical equilibrium conditions.
Here $\bm{u}(\bm{x})$ is the displacement field, $\mathbb{C}(\bm{x})$ is the fourth order elastic modulus tensor, $b(\bm{x})$ is the body force, $\bm{u}^0(\bm{x})$ is the displacement specified at the Dirichlet boundary $\partial_1\Omega$, $\bm{t}^0(\bm{x})$ is the traction specified in the traction boundary $\partial_2\Omega$, and $\hat{\bm{n}}(\bm{x})$ is the normal vector at any point on the boundary. 
We also assume the conventional major and minor symmetries of $\mathbb{C}(\bm{x})$ to allow us to directly work with the displacement $\bm{u}(\bm{x})$.
While the Dirichlet boundary condition is essential to solid mechanics problems, imposing them on the diffused boundary has limited use cases.
Therefore, we limit our attention to the traction boundary condition.
To move to the diffuse boundary setting, we introduce the diffuse traction $\hat{\bm{t}}^0:\partial_\epsilon\Omega_\epsilon\to\mathbb{R}^3$,
\begin{equation}
  \hat{\bm{t}}^0(\bm{x}=\bm{y} + s\bm{n}) = \bm{t}^0(\bm{y}).
\end{equation}
The diffuse traction is defined everywhere in the diffuse boundary as the value of the discrete traction at the closest point on the discrete boundary.
(Recall that this is only valid as long as $\epsilon$ is smaller than the smallest radius of curvature of the discrete boundary; otherwise, the diffuse traction is multiply defined.)
Now, consider the following diffuse-boundary modification of \cref{eq:momentum}:
\begin{equation}
  \big( \operatorname{div} \bm{\sigma} - \bm{b} \big)\phi^\epsilon  = \big( \hat{\bm{t}}^0 - \bm{\sigma}\bm{n} \big) |\nabla\phi^\epsilon|. \label{eq:momentum2}
\end{equation}
It is straightforward to show that the interior momentum equation holds by considering the weak form of the above equation.
Integrate both sides over an arbitrary interior, measurable region $V\subset\Omega_\epsilon$,
\begin{equation}
  \int_V \big( \operatorname{div} \bm{\sigma} - \bm{b} \big) \phi^\epsilon\, dV = 0 \ \ \ \forall \operatorname{meas.} V \subset \Omega_\epsilon\label{eq:mom_int_weak}.
\end{equation}
The right hand side vanishes since $|\nabla\phi^\epsilon|$ is zero in $\Omega_\epsilon$, by construction.
But since \cref{eq:mom_int_weak} holds for all subsets of the diffuse interior, the integrand itself is zero for all $\bm{x}\in\Omega_\epsilon$.
To show recovery of the boundary condition, we once again take the weak form of \cref{eq:momentum2}, this time over an arbitrary region within the diffuse boundary:
\begin{align}
  \int_A\int_{-\epsilon/2}^{\epsilon/2} \Big(\big( \operatorname{div} \bm{\sigma} - \bm{b} \big)\phi^\epsilon - \big( \hat{\bm{t}}^0 - \bm{\sigma}\bm{n} \big) |\nabla\phi^\epsilon|\Big)\,ds\,dA = 0 \ \ \  \forall A\subset\partial\Omega\label{eq:eq_mom_bc_weak}
\end{align}
Applying Theorem 1 in the sharp interface limit reduces the above expression to
\begin{align}
  \int_A\big(\bm{t}^0 - \bm{\sigma}\bm{n} \big)\,dA = 0 \ \ \  \forall A\subset\partial\Omega,
\end{align}
which shows that the integrand must be true for all $\bm{y}\in\partial\Omega$ since the above weak form holds for all subsets of $\partial\Omega$.
This confirms that the traction boundary conditions are exactly recovered as $\epsilon\to0$ for the diffuse interface formulation in \cref{eq:momentum2}.
Finally, rearranging and noting that $\bm{n} = \nabla\phi^\epsilon / |\nabla\phi^\epsilon|$., \cref{eq:momentum2} simplifies to
\begin{equation}
  \operatorname{div} \left( \phi^\epsilon \bm{\sigma} \right) - \phi^\epsilon \bm{b} = \hat{\bm{t}}^0 |\nabla\phi^\epsilon|. 
  \label{eq:SBMLinearElasticity}
\end{equation}
In summary, the boundary conditions associated with the discrete natural boundary $\partial\Omega$ are replaced by an equivalent source term that mimics the effect of the discrete natural boundary, exactly recovering its behavior in the sharp interface limit.
The selection of $\phi$ can be determined by construction (for instance, explicitly prescribing an indicator function based on a predetermined geometry) or by coupling to a separate set of equations that describe the behavior of $\phi$ (such as phase field).
In both cases, care must be taken that the behavior of $\phi$ does not deviate far from the requirements necessary for the validity of the diffuse boundary method to hold.

\subsection{Reflux-free multigrid implementation}\label{sec:refluxfree}
We implemented equation (\ref{eq:SBMLinearElasticity}) in an in-house code, Alamo \cite{runnels2021massively}, a finite-difference based multi-level, multi-grid, and multi-component solver.
Alamo uses AMReX  libraries for block-structured adaptive mesh refinement (BSAMR) \cite{zhang2019amrex}.
BSAMR divides the mesh into levels such that each level contains cells of the same size. 
Each level is treated independently, and the information between levels is communicated through restriction, relaxation, and restriction operations using ghost cells.
As a result, BSAMR is highly scalable and enables massive parallelism across CPU and GPU cores.
BSAMR can be naturally combined with standard multigrid methods by treating refined levels as an extension to the multigrid method's coarse/fine level sequence.
The mesh refinement and coarsening are triggered and performed at regular intervals using the Berger-Rigoutsos algorithm \cite{berger1991algorithm}.

Multigrid methods often require special treatment at the coarse-fine level interface during restriction operations.
Improper handling of the coarse-fine interface can result in spurious forces at the interface and overall poor convergence of the solver.
The coarse-fine interface can be handled using a ``reflux'' operation \cite{almgren1998conservative} where the operator is updated at the interface to use the information at both levels.
However, this process can be difficult for a complicated operator such as the one for linear elasticity.
An alternate ``reflux-free'' procedure was proposed by \cite{runnels2021massively} where the levels are padded with an extra layer of ghost nodes/cells to ensure the translational symmetry of the restriction operator and that information at the coarse/fine boundary is updated with the current information.
This circumvents the need for a special stencil at the coarse/fine boundary and results in good convergence.

\subsection{Node-based and cell-based fields in strong form multigrid elasticity}
The method discussed in \cref{sec:refluxfree} works well for elasticity problems but can behave very poorly unless care is taken to respect the proper placement of the relevant fields on the grid.
In fluid mechanics, it is often necessary to place some quantities at nodes, some in cells, and some on cell faces, etc., with the exact scheme differing between methods \cite{piller2004finite,alves2021numerical,wesseling1991finite,anderson1995computational}.
On the other hand, in solid mechanics, values such as displacements are typically stored at nodes, whereas quantities governing material response are generally located at quadrature points within the element.
Some solid mechanics methods, such as optimal transport meshfree \cite{li2010optimal}, smoothed particle hydrodynamics \cite{monaghan1992smoothed}, and the material point method \cite{sulsky1994particle,sulsky1995application,liang2019efficient}, though not strictly finite element, still carefully distinguish between nodes and material points.
The finite volume method has been used for solid mechanics, though not nearly as extensively as in fluid mechanics, and it is known that a staggered grid approach is needed to avoid the phenomenon of checkerboarding \cite{cardiff2021thirty}.

In the present method, which uses a regular cartesian grid, values may be stored at points, edges, faces, or cells.
At first glance, there is no obvious reason for storing values at one location over another; indeed, it is possible to develop a solver in which all values are stored at faces or all values at nodes.
In previous work, nodal locations were chosen for all quantities of interest \cite{runnels2021massively}.
Interestingly, this choice had no previous negative impact on the solver.
However, as we will discuss here, the location of values is, in fact, quite essential to the performance of the solver when considering near-singular problems.
Specifically, it is important that displacements be stored at the nodes, whereas the order parameter must be stored as a cell-based field.
In this section, we provide two explanations for this: the first, from a practical perspective, and the second, by considering the geometry of the problem.

Our analysis considers the key aspect of the multigrid solver as applied to the elasticity problem: smoothing.
Multigrid methods work primarily through the use of a smoother between levels.
The choice of smoothing algorithm can vary, but the preeminent solvers are generally Gauss-Seidel, Jacobi, or a variant of one of these methods.
These methods are popular due to their ease of implementation, their parallel efficiency.
For geometric multigrid methods, they are particularly attractive because they smooth high frequency error faster than low frequency error (unlike, for instance, the conjugate gradient method) \cite{wesseling1995introduction}.

Here, we consider the Jacobi method, which for an operator $A$ is given by
\begin{equation}
  \bm{u}^{n+1} = D^{-1}(\bm{b} - (A-D)\bm{u}^{n}),
\end{equation}
where $\bm{u}^n$ is the solution at iteration $n$, $D=\operatorname{diag}(A)$, and $\bm{b}$ is the right hand side.
We see that the inverse of the diagonal is of central importance: most notably, that if any of the diagonal elements of $D$ are zero, the corresponding rows of $(\bm{b}-(A-D)\bm{u}^n$) must be zero as well.
To see the significance of this, consider the discretized elastic operator for a constant modulus $C$ in one dimension, in which all values are stored at nodes:
\begin{align}
  \operatorname{div} \left(\phi^\epsilon \mathbb{C} \operatorname{grad} \bm{u} \right)_i &= C \left(\frac{\phi^\epsilon_{i+1} - \phi^\epsilon_{i-1}}{2\Delta x} \right) \left(\frac{u_{i+1}-u_{i-1}}{2\Delta x}\right) + C \phi^\epsilon_i \left(\frac{u_{i+1} - 2u_i + u_{i-1}}{\Delta x^2}\right) 
\end{align}
The diagonal of the operator, then, is nothing other than the coefficient of the $u_i$ term, that is,
\begin{align}
  \operatorname{diag}(\operatorname{div} [\phi_\epsilon \mathbb{C} \operatorname{grad}] )_i
  &=  -2\frac{C \phi^\epsilon_i}{\Delta x^2},
\end{align}
and so the corresponding Jacobi update is thus given by
\begin{align}
  u_i^{n+1} 
  &= -2\frac{\Delta x^2}{C \phi^\epsilon_i}
    \Bigg\{
    b_i - 
    \Big[ C \Big(\frac{\phi^\epsilon_{i+1} - \phi^\epsilon_{i-1}}{2\Delta x} \Big) \Big(\frac{u^n_{i+1}-u^n_{i-1}}{2\Delta x}\Big) + C \phi^\epsilon_i \Big(\frac{u^n_{i+1} + u^n_{i-1}}{\Delta x^2}\Big) \Big]
    \Bigg\}\\
  &= -2\Delta x^2
    \Bigg\{
    \frac{b_i}{\phi_i^\epsilon C} - 
    \Big[ \frac{1}{\phi^\epsilon_i}\Big(\frac{\phi^\epsilon_{i+1} - \phi^\epsilon_{i-1}}{2\Delta x} \Big) \Big(\frac{u^n_{i+1}-u^n_{i-1}}{2\Delta x}\Big) + \Big(\frac{u^n_{i+1} + u^n_{i-1}}{\Delta x^2}\Big) \Big]
    \Bigg\}
\end{align}
By inspection, it is clear that any nonzero body force will produce divergent behavior if applied where $\phi_i^\epsilon=0$; this is natural since such a problem would be ill-defined.
However, an inspection of the next term shows a second vulnerability: a point at which $\phi^\epsilon_i=0$ may still induce instability depending on the values at {\it the adjacent nodes}.
Therefore this can (and does) induce divergence at the nodes where the order parameter is zero, but the solution is well-defined: i.e., at the boundaries of the support of $\phi^\epsilon$.

On the other hand, consider the corresponding Jacobi update if the field $\phi^\epsilon$ is stored in cells rather than at nodes, where fractional indices are used to denote locations of cells:
\begin{align}
  u_i^{n+1} 
  = -2\Delta x^2
  \Bigg\{&
           \frac{2b_i}{(\phi_{i+1/2}^\epsilon +\phi^{\epsilon}_{i-1/2})C} \notag\\
         &- 
           \Big[ \frac{2}{(\phi_{i+1/2}^\epsilon +\phi^{\epsilon}_{i-1/2})}\Big(\frac{\phi^\epsilon_{i+1/2} - \phi^\epsilon_{i-1/2}}{\Delta x} \Big) \Big(\frac{u^n_{i+1}-u^n_{i-1}}{2\Delta x}\Big) + \Big(\frac{u^n_{i+1} + u^n_{i-1}}{\Delta x^2}\Big) \Big]
           \Bigg\}.
\end{align}
One can see by inspecting the second term (assuming, again, that the body force is responsibly applied) that the problem of instability is eliminated.
Divergent behavior can now only occur when $\phi^{\epsilon}$ is zero at both $i+1/2$ and $i-1/2$; but if this happens, then the difference between the two values would be zero as well.
One may apply this same exercise to this problem with non-uniform elastic modulus, or to the problem in 2D or 3D, with the same result.
This underscores the importance of a staggered grid approach that, though commonly used in other finite difference methods, was absent from prior finite difference implementations of the SBM on BSAMR grids.

One may take this analysis several steps further by considering the geometric significance of the solution and $\phi^\epsilon$ fields.
Within the past couple of decades, the tools of exterior calculus have been applied to the problem of linear elasticity \cite{kanso2007geometric,yavari2008geometric,yavari2010geometric}, which identifies displacement fields as vector-valued 1-forms, body forces as vector-valued 3-forms, etc.
This has been extended to the field of computational mechanics through the emerging sub-discipline of discrete differential geometry (DDG), which allows the explicit realization of exterior calculus constructs in the context of discrete mesh elements \cite{desbrun2006discrete,zhu2022mem3dg,ruocco2021discrete}.
While a thorough treatment of DDG in the context of mechanics is outside the scope of the paper, we outline the underlying concept.
The traditional fields (displacements, stress, etc) can be replaced by the DDG construct of forms, where an n-form is a field that can be integrated over an n-dimensional manifold.
That is, n-forms contain a notion of geometry that is absent from raw fields, and which correspond to the proper integration domain.
For instance, stress are integrated over surfaces, which are two-dimensional manifolds; therefore, stress is a 2-form;
body forces are integrated over volumes, which are three-dimensional manifolds; therefore, body forces are 3-forms;
displacements are evaluated at points, which are zero-dimensional manifolds; therefore, displacements are 0-forms; and so on.
In the present work, the order parameter $\phi$ is integrated over a voluem; therefore, $\phi$ is a 3-form.
In the discrete setting on a regular grid, 0-forms may be identified as nodal fields, 1-forms as edge fields, 2-forms as face fields, and 3-forms as cell fields.
(The study of DDG has produced analagous interpretations for non-regular and unstructured meshes as well.)
Thus, we see that by representing displacements using a node-based grid (0-forms) and $\phi$ using a cell-based grid (3-forms), we are preserving the geometric structure of the these fields.
For a more thorough discussion of the geometric interpretation of classical continuum theory, we refer the reader to the above references and to the excellent book by Frenkel \cite{frankel2011geometry}.

\subsection{Material model vector space}
One of the key advantages of the BSAMR approach is its ability to adapt the mesh rapidly.
The block-structured multilevel data structures afford rapid regridding in a parallel-efficient manner.
Regridding, and inter-level communication, relies on the ability to rapidly transfer information between AMR (or multigrid) levels, usually in the form of interpolation or prolongation.
When working with primitive fields such as velocity, density, pressure, etc., interpolations and prolongations are easy to compute.
However, such operations are not always obvious in the context of material modeling.
Materials often exhibit highly anisotropic behavior, with material response often depending on numerous parameters and time-evolving internal variables.
In order for BSAMR to function correctly with such models as these, it is necessary to address the constraints and requirements needed for material modeling.
Moreover, as it is a commonplace in any solid mechanics code to allow for modular material models, an {\it ad hoc} implementation is insufficient.
Therefore, we prescribe the minimum requirements needed for a versatile implementation of material models in a BSAMR context.

The aforementioned requirements for BSAMR data structures are equivalent to those for the mathematical structure of a vector space.
Specifically, BSAMR requires the consistency of solid models between AMR and multigrid levels, which is achieved through interpolation and restriction operations, which require the definition of addition and scalar multiplication of solid model objects.
Therefore we impose the requirement on material models that they must satisfy the properties of a vector space.
Let $\mathcal{M}$ denote the vector space corresponding to a certain material model.
The salient properties are: (1) the existence of a ``vector addition'' operation, typically denoted $+$, such that $a+b\in\mathcal{M}$ $\forall a,b\in\mathcal{M}$; and (2) the existence of a ``scalar multiplication'' operation denoted by ``*'' or concatenation, such that $\alpha * a = \alpha a \in \mathcal{M}$ $\forall \alpha\in\mathbb{R},\forall a \in \mathcal{M}$, and (3) the existence of an identity element $0\in\mathcal{M}$ such that $0+a=a$ $\forall a \in \mathcal{M}$.
Other requirements include associativity and commutativity of $+$, the inverse of $+$, compatibility of $*$ and identity under $*$, and distributivity of $*$; generally, these are not troublesome to enforce, and they furnish a valuable framework for unit testing at the implementation phase.
Properties 1-3 must hold not only for the internal variables stored in each material model but for their functions as well: specifically, the zeroth, first, and second derivatives of energy (W, DW, DDW), as well as any functions defining the evolution of internal variables.
For instance, $(a+b).DW(\mathbf{\varepsilon}) \overset{!}{=} a.DW(\bm\varepsilon) + b.DW(\bm\varepsilon)$, where $\bm\varepsilon$ is the local strain tensor.
The structure also allows for the inverse of models to exist, which allows for derivatives of models to be calculated, e.g. $(da/dx).W(\varepsilon)$.
For instance, the calculation of the gradient in the $x_1$ direction of a model field $a(\bm{x})$ using finite difference would be
\begin{align}
  \Big(\frac{da}{dx_1}\Big).DW(\bm{\varepsilon}) \approx \Big(\frac{a(x_1+\frac{1}{2}\Delta x_1,x_2,x_3)-a(x_1-\frac{1}{2}\Delta x_1,x_2,x_3)}{\Delta x_1}\Big).DW(\bm{\varepsilon}),
\end{align}
which makes use of the full agebraic structure of $a$: namely, algebraic operations of scalar multiplication and addition, as well as inversion.

The model $-a$ is the inverse of $a$ with negative (and consequently unphysical) material properties.
We emphasize the importance of placing checks in place to ensure that unphysical models are not accidentally used to calculate real physical properties.

The vector space material model requirement has immediate implications on the implementation of material models.
Consider the simple case of linear elastic isotropic, which is generally parameterized by two properties, often chosen as Young's modulus $E$ and Poisson's ratio $\nu$.
One can construct a material model based on these two properties $(E,\nu)$ along with the addition operation $(E_1,\nu_1)+(E_2,\nu_2) = (E_1+E_2,\nu_1+\nu_2)$, and scalar multiplication $\alpha(E,\nu)=(\alpha E,\alpha \nu)$.
However, such a model violates the vector space behavior of $W$, $DW$, and $DDW$, since the energy, stress, and strain depend on the ratio $\nu/E$ rather than bilinearly on $\nu$ and $E$ separately.
Thus, one can instead store the Lam\'e constants $\lambda,\mu$, on which the dependence of $W, DW, DDW$ is bilinear.

Another salient example is the implementation of cubic elasticity, which requires the storage of rotational information along with elastic moduli $C_{11},C_{12},C_{44}$.
One may also include an eigenstrain $\bm{\varepsilon}_0$, reflecting plastic evolution, thermal expansion, etc.
Euler angles are sometimes used to store the local rotation but are clearly a poor choice here, as Euler angles do not form a vector space.
Instead, we used quaternions to store rotation information, as they possess an algebraic structure that is relatively easy to implement.
One complication is that quaternions must be normalized to obtain rotation information, meaning that the $W,DW,DDW$ functions for the zero element are ill-defined.
However, in practice, it is generally never the case that the zero element would be called upon to return those values, and if it did, it would always return zero anyway by necessity.
The remaining values in the model readily admit an algebraic structure and are easily implemented.

Operator overloading (that is, the definition of functions using the syntax of operators such as +, -, +=, etc.)
was used to provide both unary and binary vector addition and scalar multiplication operations, and unit tests can enforce associativity, commutativity, etc.
We used forced code inlining (compile-time restructuring to place called code ``inline'' with the calling code) to ensure that there is no function call overhead, and C preprocessor macros can be used to implement all operators with minimal boilerplate code required automatically.
We required functions such as \verb-Zero- to furnish the zero element, with template metaprogramming used to enforce that all models comply with the vector space requirement.
We refer the interested reader to \cite{aigner1996eliminating,eijkhout2022science,abrahams2004c++} for further discussion of relevant high performance computational techniques.

\section{Examples}\label{sec:examples}
In this section, we demonstrate the performance and accuracy of the solver and the SBM implementation within Alamo using problems within solid mechanics.
\subsection{Plastic deformation due to spherical void}
\begin{figure}\centering
  \begin{subfigure}{0.5\linewidth}
    \centering
    \includegraphics[height=4.5cm]{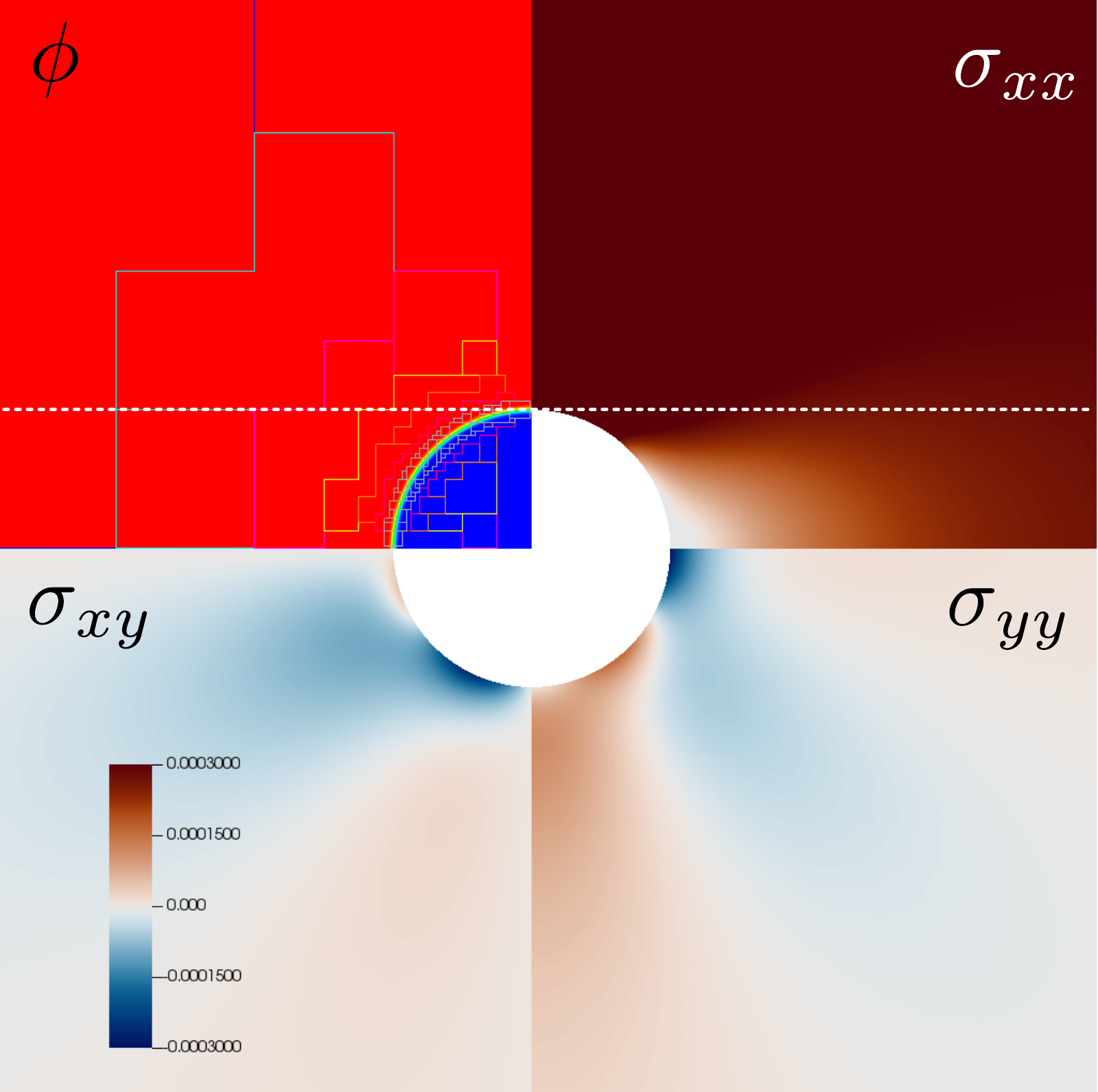}
    \caption{Plot of $\phi$ (top left) and stress fields $\sigma_{xx}$ (top right), $\sigma_{yy}$ (bottom right), $\sigma_{xy}$ (bottom left) for uniaxial tension in the x direction.}
    \label{fig:platehole_plot}
  \end{subfigure}%
  \begin{subfigure}{0.5\linewidth}
    \centering
    \includegraphics[height=4.5cm]{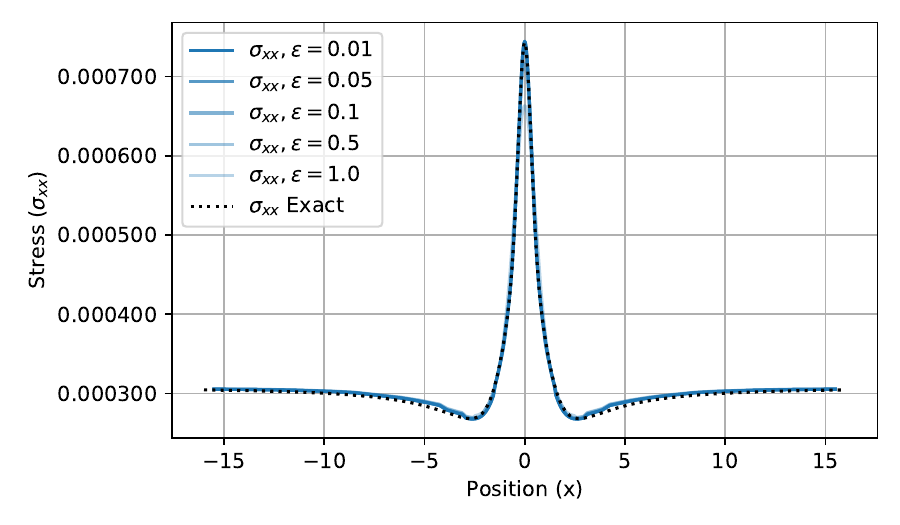}
    \caption{
      $\sigma_{xx}$
    }
    \label{fig:platehole_comparison_11}
  \end{subfigure}
  \begin{subfigure}{0.5\linewidth}
    \centering
    \includegraphics[height=4.5cm]{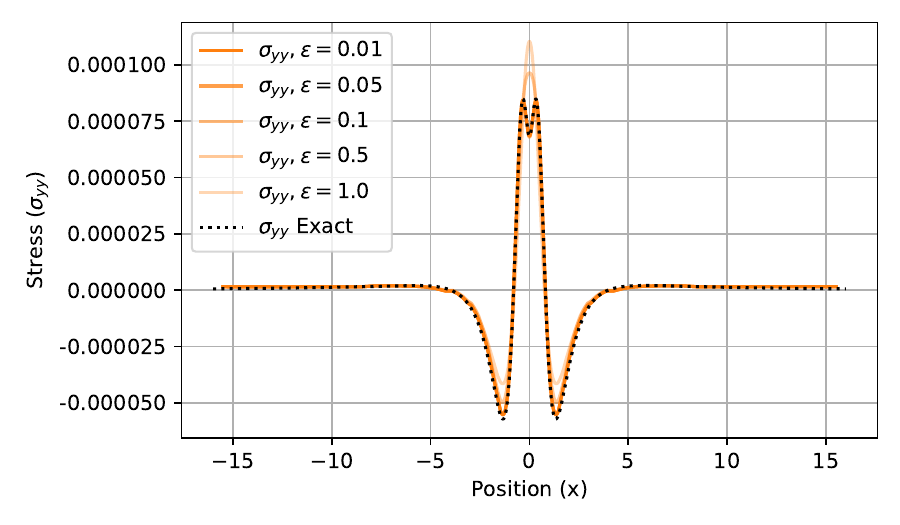}
    \caption{
      $\sigma_{yy}$
    }
    \label{fig:platehole_comparison_22}
  \end{subfigure}%
  \begin{subfigure}{0.5\linewidth}
    \centering
    \includegraphics[height=4.5cm]{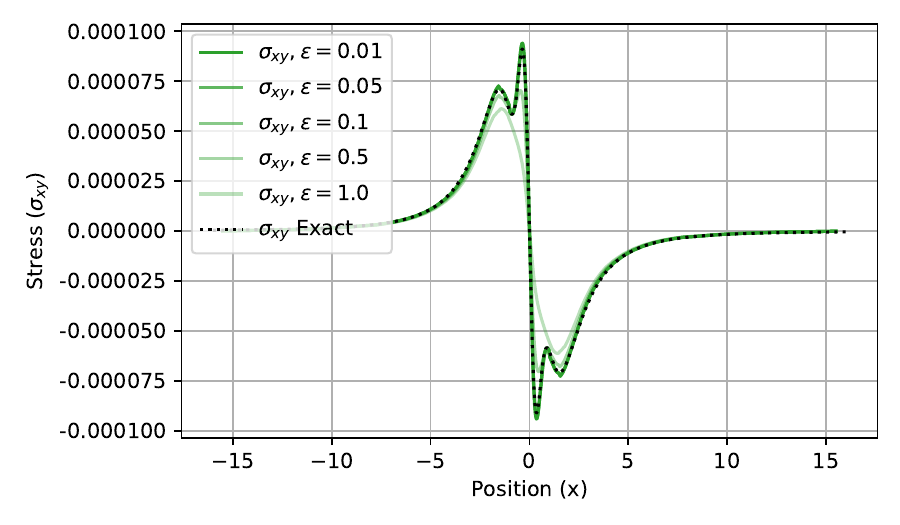}
    \caption{
      $\sigma_{xy}$
    }
    \label{fig:platehole_comparison_12}
  \end{subfigure}
  \begin{subfigure}{0.5\linewidth}
    \centering
    \includegraphics[height=4.5cm]{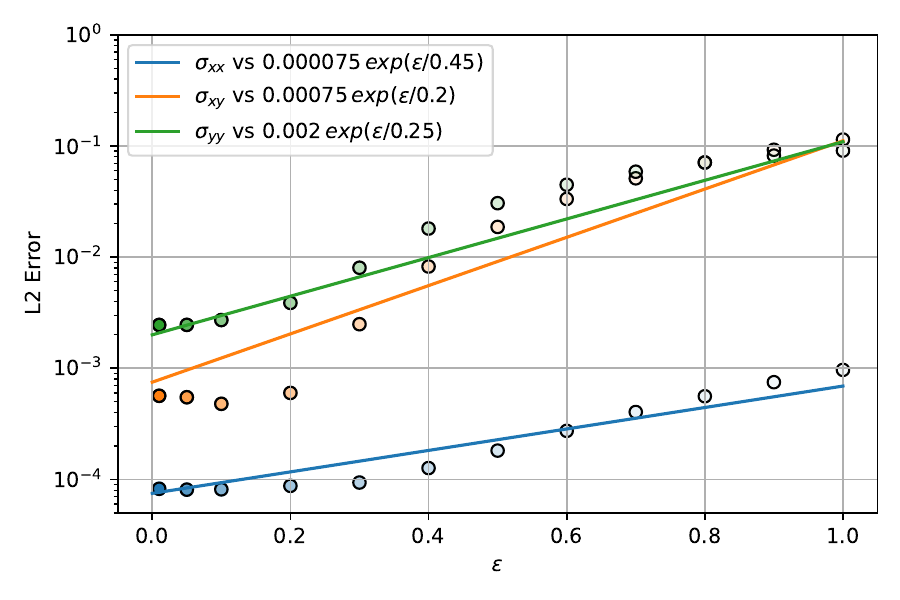}
    \caption{
      L2 error (as compared to the exact solution) with respect to diffuse boundary thickness $\varepsilon$.
    }
    \label{fig:platehole_comparison_error}
  \end{subfigure}

  \caption{
    Verification of the method by comparison to the exact solution for a plate with a hole under uniaxial tension.
    Convergence with decreasing $\varepsilon$ (corresponding to increasing decreasing line opacity), compared to the exact solution (dashed lines), along $y=R=1$, which is the tangent to the hole (dashed white line).
  }
  \label{fig:platehole}
\end{figure}

As a first step, we solved a standard canonical problem to validate the accuracy of the near singular solver.
We considered the two-dimensional problem of a large linear elastic plate with a circular hole subjected to uniaxial stress. 
The stress fields around the hole are well-defined and can be analytically computed using the Airy stress function approach \cite{bower2009applied}.
We chose a two-dimensional domain of $\bm{x}\in [-16,16]\times [16,16]$ and introduced a circular hole of radius $1.0$ at the center.
We used a range of regularization length scales  $0.01$, $0.05$, $0.1$, $0.5$, and $1.0$.
We subject the domain to a uniaxial stress condition by fixing the left edge and applying a displacement in the x direction on the right edge. 
\Cref{fig:platehole_plot} (top left) shows the $\phi^\epsilon$ field with a regularization length scale of $0.01$ along with the refined grid.
The corresponding stress distributions $\sigma_{xx}$, $\sigma_{xy}$, and $\sigma_{yy}$ are shown in \Cref{fig:platehole_plot} in the top right, bottom left and bottom left respectively.
We present the comparison of the numerical solution with the analytical solution at $y=1$ line (tangent to the hole) as a function of the regularization scale in \Cref{fig:platehole_comparison_11,fig:platehole_comparison_12,fig:platehole_comparison_22}.
We note that the solution predicted by the near-singular solver converges to the analytical solution \cite{bower2009applied} as the regularization length scale decreases.
The normalized L2 error (i.e. norm of the difference divided by the norm of the reference) indicates convergence with respect to $\varepsilon$ (\Cref{fig:platehole_comparison_error}).
The increased error for the smallest 1-2 values of $\varepsilon$ indicates an error from the discretization since a constant mesh resolution (that is, the same number of AMR levels) was used for each case for consistency.

\begin{figure}
  \centering
  \includegraphics[width=0.6\linewidth]{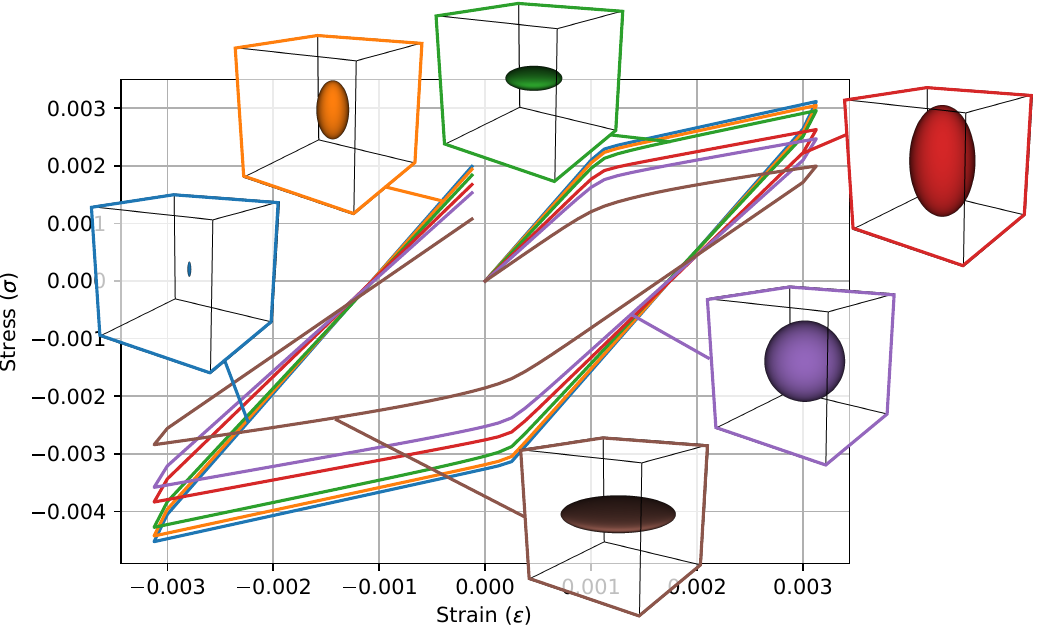}
  \caption{Stress-strain hysteresis results for a variety of void shapes and sizes.}
  \label{fig:j2_hysteresis}
\end{figure}

\begin{figure}
  \centering
  \includegraphics[width=0.33\linewidth]{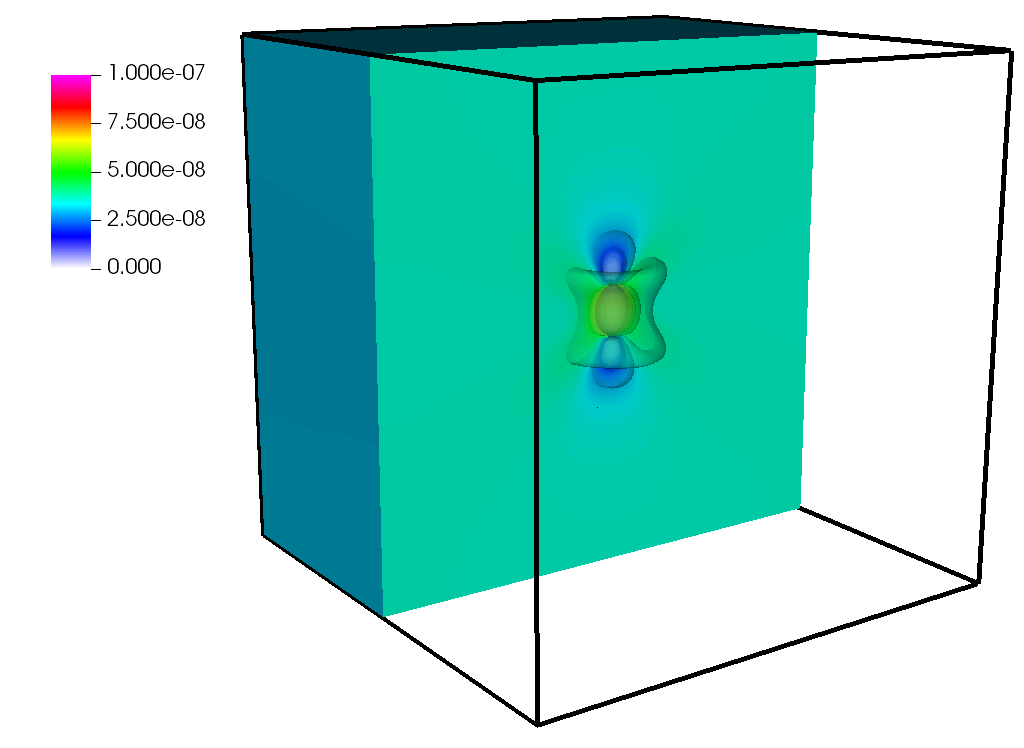}%
  \includegraphics[width=0.33\linewidth]{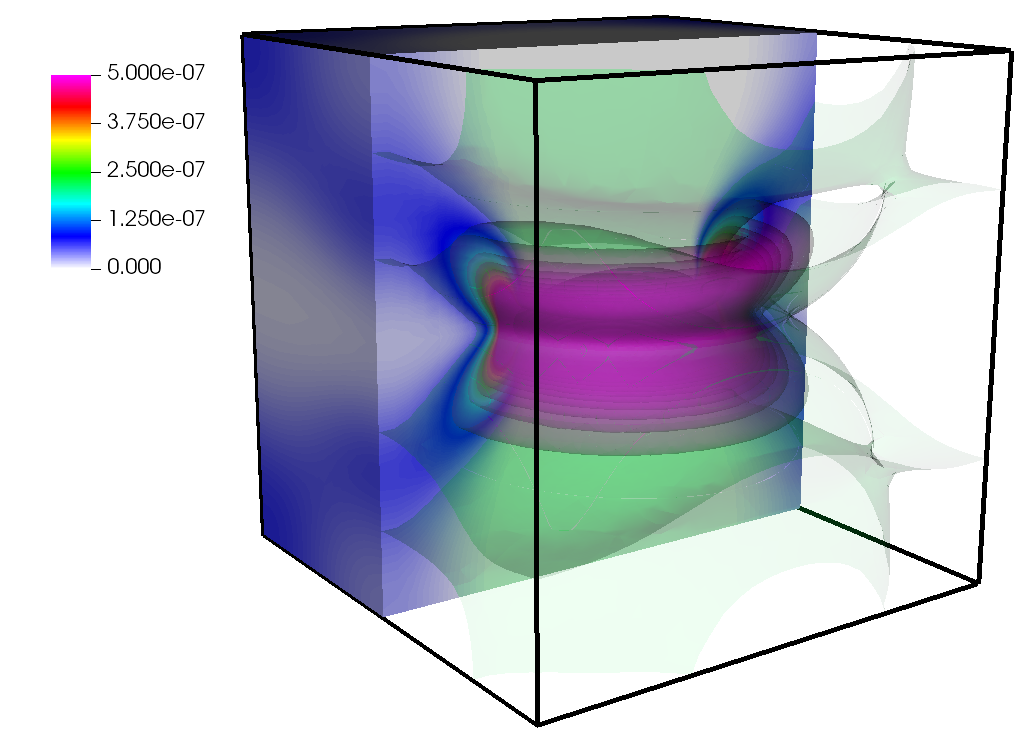}%
  \includegraphics[width=0.33\linewidth]{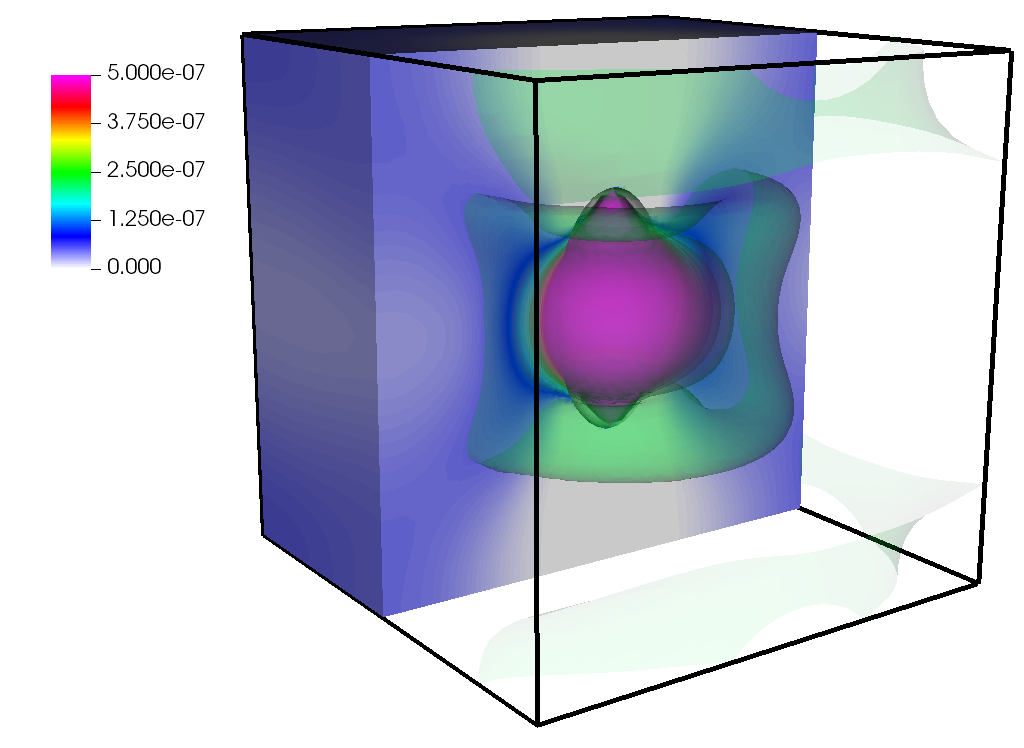}
  \includegraphics[width=0.33\linewidth]{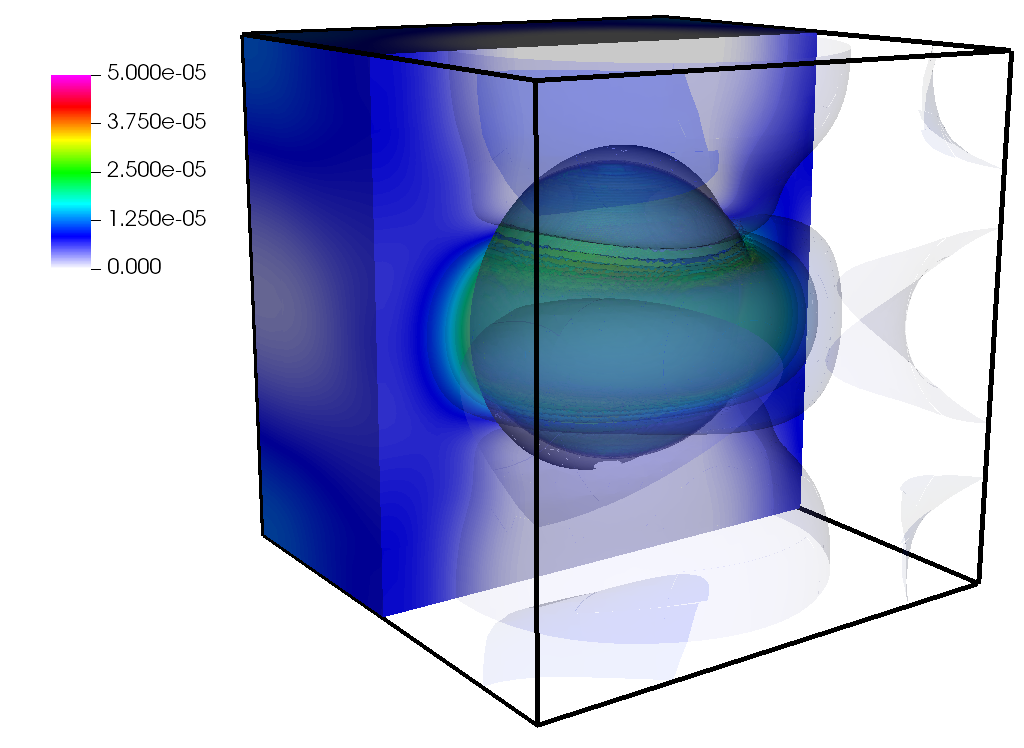}%
  \includegraphics[width=0.33\linewidth]{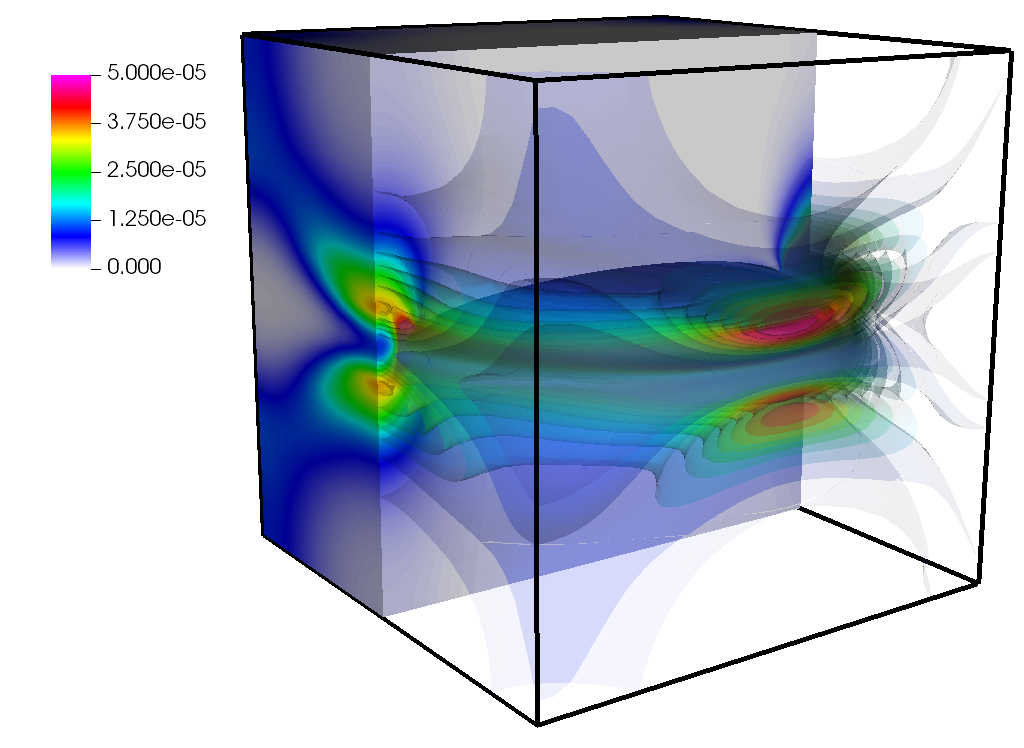}%
  \includegraphics[width=0.33\linewidth]{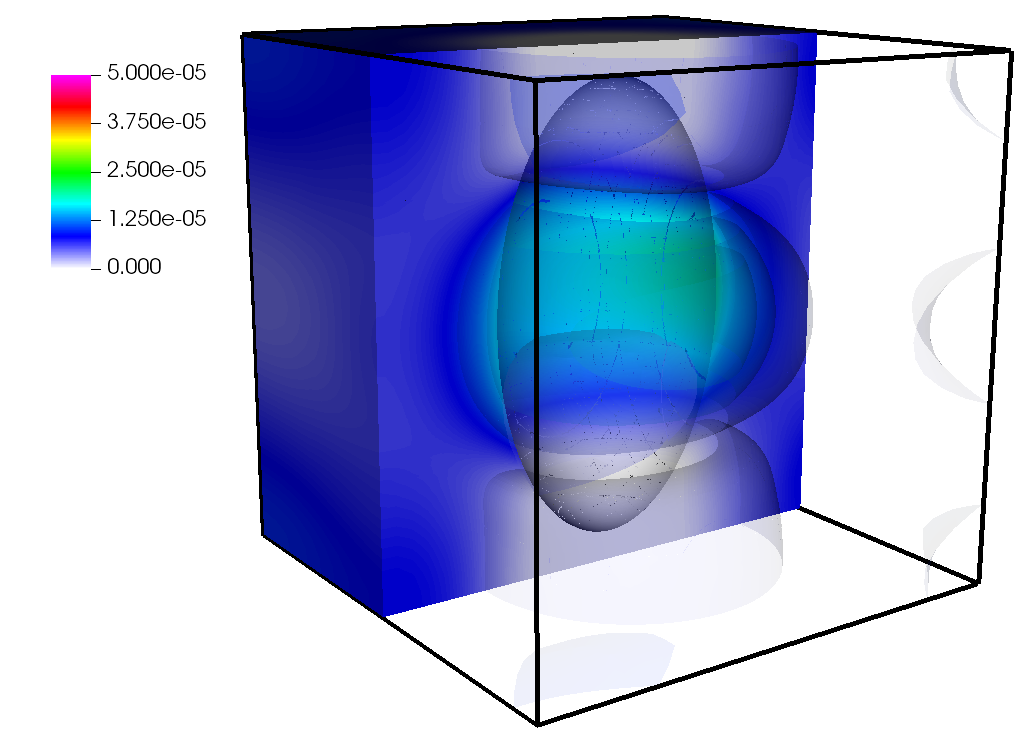}
  \caption{Isocontours showing the magnitude of the plastic strain deviator at $t=0.3$. From top right to bottom left, the void radii are $\bm{r}_1=(2.0,0.5,0.5)$, $\bm{r}_2=(2.0,7.07,7.07)$, $\bm{r}_3=(7.5,4.08,4.08)$, $\bm{r}_4=(10,10,10)$, $\bm{r}_5=(2.5,14.14,14.14)$, and $\bm{r}_6=(14.14,8.16,8.16)$.}
  \label{fig:J2_epspvisualization}
\end{figure}

Having validated the elastic solver, we demonstrate the effectiveness of the method by considering the plastic deformation of a cuboidal object with an embedded spherical void subject to uniaxial loading.
We represented the material using the order parameter $\phi_\epsilon$ which takes the value $0$ within the void and $1$ outside, with a length scale $\epsilon$.
The stress in the SBM equation (\ref{eq:SBMLinearElasticity}) is expressed as $\sigma = \mathbb{C}\left( \operatorname{grad} \bm{u} - \bm{\varepsilon}_p \right)$, where $\bm{\varepsilon}_p$ is the plastic strain.
We used a staggered approach to solve the elastic equilibrium equation (\ref{eq:SBMLinearElasticity}) and plastic evolution, and modeled the evolution of plastic strain $\bm{\varepsilon}_p$ using the $J_2$-plasticity model for a linear elastic isotropic material with Lam\'e constants $\lambda$ and $\mu$.
We chose a $J_2$ plastic strength model with isotropic hardening.
The yield strength for isotropic hardening is given by
\begin{equation}
  K(\alpha) = \sigma_Y + \theta \bar{H}\alpha, \quad \theta\in[0,1]
\end{equation}
where $\sigma_Y$ is the flow stress, $\alpha$ is the equivalent plastic strain, $\bar{H}$ is the hardening modulus and $\theta$ is a parameter governing the hardening slope.
We used the internal variables $\bm{q} = \{\alpha, \bm{\varepsilon}_p, \bm{\beta}\}$ for the plasticity model, where $\beta$ is the center of the von-Mises yield surface in the stress deviator space.
Following are the yield condition flow rule and hardening rule for the $J_2$ plasticity model. 
\begin{align}
  \bm{\eta} &:= \operatorname{dev}[\bm{\sigma}] - {\bm{\beta}}, \quad \operatorname{tr} {\bm{\beta}} := 0,\quad 
              f(\bm{\sigma},\bm{q}) = ||\bm{\eta} || - \sqrt{\frac{2}{3}} K(\alpha), \nonumber \\
  \dot{\bm{\varepsilon}}_p &= \gamma \frac{\bm{\eta}}{||\bm{\eta}||}, \quad \dot{\alpha} = \gamma\sqrt{\frac{2}{3}}, \quad  \dot{\bm{\beta}} = \gamma\frac{2}{3}(1-\theta) \bar{H} \frac{\bm{\eta}}{||\bm{\eta}||} \label{eq:J2Plasticity}
\end{align}
We solved the above equations using the following radial return algorithm described in detail in \cite{simo2006computational}.
Given a stress and strain state at time $t_n$ as $\bm{\sigma}_n$ and $\bm{\varepsilon}_n$, and the strain $\bm{\varepsilon}_{n+1}$ at time $t_{n+1}$, the algorithm involves following steps.

\begin{algorithm}[H]
  \caption{J2 plasticity model update}
  \begin{algorithmic}[1]
    \State 
    $\bm{e}_n \gets \operatorname{dev}\bm{\varepsilon}_n$, 
    $\bm{e}_{n+1} \gets \operatorname{dev}\bm{\varepsilon}_{n+1}$, 
    $\bm{s}_{n} \gets \operatorname{dev}\bm{\sigma}_{n}$, 
    \Comment{Compute deviatoric strain and stress}
    \State
    $\bm{s}_{n+1}^\mathrm{trial} \gets \bm{s}_n + 2\mu(\bm{e}_{n+1}-\bm{e}_n)$,
    $\bm{\eta}^\mathrm{trial}_{n+1} \gets \bm{s}^\mathrm{trial}_{n+1}-\bm{\eta}_n$
    \Comment{Compute trial states}
    \State
    $\bm{n}_{n+1} \gets \bm{\eta}^\mathrm{trial}_{n+1}/|\bm{\eta}^\mathrm{trial}_{n+1}|$
    \Comment{Compute new yield surface normal}
    \State
    Solve $-\sqrt{\frac{2}{3}}K\left(\alpha_n + \sqrt{\frac{2}{3}} \Delta \gamma\right) + ||\bm{\eta}_{n+1}^{trial} || = 0$ for $\Delta\gamma$
    \Comment{Consistency condition}
    \State 
    $\bm{\varepsilon}_{n+1}^p = \bm{\varepsilon}_n^p + \Delta\gamma\bm{n}_{n+1} $
    \Statex
    $\alpha_{n+1} = \alpha_n + \sqrt{\frac{2}{3}} \Delta \gamma$
    \Statex
    $\bm{\beta}_{n+1}^p = \bm{\beta}_n^p + \sqrt{\frac{2}{3}}\,\theta\, \bar{H} (\alpha_{n+1}-\alpha_n)\bm{n}_{n+1}$
    \Comment{Update internal variables}
  \end{algorithmic}
\end{algorithm}

A 3D domain $\bm{x}=(x_1,x_2,x_3)\in[-16,16]\times[-16,16]\times[-16,16]$ (arbitrary units) with an ellipsoidal void centered at the origin and radii $\bm{r} = (r_x, r_y, r_z)$ was used. 
The order parameter $\phi$ was set to $1$ outside the inclusion and $0$ inside with a length scale $\epsilon=0.4$.
We chose the material parameters as Young's modulus $E=210\,\text{GPa}$, Poisson's ratio $\nu=0.3$, yield strength $\sigma_Y=200\,\text{MPa}$, and hardening parameters $\bar{H}=50\,\text{GPa}$ and $\theta=1$.
We performed a tension test with a fixed $x_1=-16$ face and a cyclic displacement applied in the $x_1$ direction on the $x_1=16$ face.
We chose the applied displacements in the increments of $0.004$ going from total applied displacement from $0.0$ to $0.1$, then from $0.1$ to $-0.1$, and finally from $-0.1$ to $0.0$.

\Cref{fig:j2_hysteresis} shows the stress-strain curves obtained for six different ellipsoid shapes and sizes.
These are $\bm{r}_1=(2.0,0.5,0.5)$, $\bm{r}_2=(2.0,7.07,7.07)$, $\bm{r}_3=(7.5,4.08,4.08)$, $\bm{r}_4=(10,10,10)$, $\bm{r}_5=(2.5,14.14,14.14)$, and $\bm{r}_6=(14.14,8.16,8.16)$.
We calculated the strains using the applied displacement and stresses from the total traction on the $x=16$ face where the displacement is applied.
As expected, the stress-strain curve for each case exhibits the classic hysteresis loop. 
As the size and aspect ratio of the void change, the plastic evolution within the domain changes leading to different stress-strain curves. 
The total plastic strain is higher for larger void with higher aspect ratios aligned with the loading directions, making the stress-strain curve flatter.
\Cref{fig:J2_epspvisualization} shows the magnitude of plastic strain deviator at the applied displacement of $0.08$ during the unloading cycle. 

We performed these simulations on the UCCS INCLINE cluster using 128 cores on a single node. 
We chose a base mesh of $32\times 32 \times 32$ with 5 levels of refinement.
For the six cases presented, the solver took $4$ minutes to $8$ hours, depending on the size of the inclusion and the size of the portion of the domain refined with high resolution. 
The solver converged linearly for all cases despite large regions of voids within the domain.
Therefore the solver performed well in predicting the stress fields and plastic strains due to voids.

\subsection{Brittle fracture}
Fracture is one of the most prominent causes of failure for engineering structures.
As such computational modeling of crack nucleation and propagation in engineering materials is critical for evaluating their performance.
Computational methods for modeling fracture can be broadly classified as either discrete boundary or diffused boundary methods.
Among the discrete boundary methods, there are two main approaches: the eXtended Finite Element Method (XFEM) \cite{moes1999finite,sukumar2000extended,li2018review} and the Scaled Boundary Finite Element Method (SBFEM) \cite{song1997scaled,song2018review,wolf2000scaled,song2000scaled}.
XFEM involves enriching classical finite elements with specialty elements designed specifically for capturing singularities at crack tips.
On the other hand, SBFEM uses a dimensional reduction technique to reduce the problem domain to the boundary of the solid and scales the solution to the crack tip analytically.
A detailed review of discrete methods can be found here \cite{egger2019discrete,sedmak2018computational}.
While these methods have been widely successful, they suffer from limitations when explicitly tracking crack fronts for complex crack patterns.

Diffuse boundary methods, or “phase field methods” use a smoothly varying scalar damage field $c(\bm{x},t)$ to diffuse the sharp crack over a length scale $\xi$ \cite{ambati2015review}.
The differential equation governing the evolution of $c(\bm{x},t)$ is based on a rigorous variational approach to fracture which uses an energy functional regularized over the length scale $\xi$ \cite{francfort1998revisiting}.
The variational method has been shown to be consistent with linear elastic fracture mechanics under quasi-static loading \cite{francfort1998revisiting} for brittle materials. 
Recently, phase field fracture methods have been extended to study heterogeneous materials \cite{nguyen2016phase, hansen2020phase, agrawal2021block}, anisotropic materials \cite{teichtmeister2017phase, li2015phase}, functionally graded materials \cite{doan2016hybrid, dinachandra2020phase, kumar2021phase}, dynamic fracture \cite{karma2001phase,borden2012phase,hofacker2012continuum}, ductile fracture \cite{ambati2015phase, miehe2015phase, kuhn2016phase}, and fatigue loading \cite{mesgarnejad2019phase}. 
Phase field fracture has also been used to study interfacial strength in composites by incorporating cohesive zone elements into the formulation \cite{carollo2018modeling, tarafder2020finite, quintanas2020phase}.
Modifications to the traditional formulation have also been studied with different damage degradation functions and damage energy penalty terms \cite{pham2011issues, kuhn2015degradation}. 

While phase field methods have gained wide adoption, they suffer from high computational costs due to typically small values of $\xi$.
One way to circumvent this challenge is to use spectral methods, typically the Fast Fourier Transform (FFT), for memory efficiency \cite{chen2019fft,ernesti2020fast}.
Another more widely used approach is to couple the phase field fracture implementation to AMR with higher resolutions near the crack tip.
Among the AMR approaches, the discontinuous Galerkin approach has been used for single-level AMR \cite{muixi2020hybridizable} and the finite cell method together with h and p refinement has been used to achieve multiple levels of refinement on a regular grid \cite{nagaraja2019phase}. 
Other attempts include hybridizing phase field method with XFEM \cite{giovanardi2017hybrid}.

In this work, we implemented a hybrid model of phase field brittle fracture to study crack propagation in Mode-I loading.
We use a regularized field $c(\bm{x},t)$ with values $1$ outside the crack and $0$ inside the crack and length scale parameter $\xi$. 
The phase field fracture energy functional is given by,
\begin{equation}
  \mathcal{L}= \int_\Omega \left(g(c) + \eta \right)W_0(\bm{\varepsilon}(\bm{u})) dV + \int_\Omega G_c \left[ \frac{w(c)}{4\xi} + \xi |\nabla c|^2\right] dV,
\end{equation}
where $W_0$ is the elastic strain energy of the material without a crack field, $G_c$ is the fracture energy, and $\eta=10^{-4}$ chosen for computational stability. 
The interpolation function $g(c)$ takes the value $0$ inside the crack and $1$ outside. 
The interpolation function $w(c)$ takes the value $1$ inside the crack and $0$ outside.
In this work, we choose a quartic interpolation function \cite{kuhn2015degradation} with $g(c) = 4c^3-3c^4$ and $w(c) = 1-g(c)$.
The crack and displacement fields are evolved using the variational derivative of $\mathcal{L}$.

We chose a linear elastic isotropic material with Lam\'e constants $\lambda$ and $\mu$ with the strain energy density and stress as
\begin{equation}
  W_0 = \frac{1}{2} \lambda (\operatorname{tr}\bm{\varepsilon})^2 + \mu \operatorname{tr}\left(\bm{\varepsilon}^2\right), \quad \bm{\sigma} = \frac{\partial W_0}{\partial \bm{\varepsilon}}
\end{equation}
To account for the tension-compression asymmetry, we assume an additive decomposition of the strain energy $W_0 = W_0^+ + W_0^-$ which uses the spectral decomposition of the strain tensor $\bm{\varepsilon} = \sum_{i=1}^d \varepsilon_i \hat{\bm{v}}_i\otimes\hat{\bm{v}}_i$.
The strain energies are given by
\begin{equation}
  W_0^\pm = \frac{1}{2} \lambda (\operatorname{tr}\bm{\varepsilon}_\pm)^2 + \mu \operatorname{tr}\left(\bm{\varepsilon}_\pm^2\right), \quad \bm{\varepsilon_\pm} = \sum_{i=1}^d \left(\varepsilon_i \right)_\pm \hat{\bm{v}}_i\otimes\hat{\bm{v}}_i,
\end{equation}
and the energy functional is updated to 
\begin{equation}
  \mathcal{L}= \int_\Omega \left[\left(g(c) + \eta \right)W_0^+ + W_0^-\right] dV + \int_\Omega G_c \left[ \frac{w(c)}{4\xi} + \xi |\nabla c|^2\right] dV.
\end{equation}
We note that the above strain-based decomposition of energy does not always ensure that the compressive states do not contribute to crack growth \cite{lo2019phase,sun2021poro} and a stress based approach has been used recently to circumvent this issue \cite{clayton2022stress}.
The equations for the hybrid formulation of phase-field fracture are obtained by taking the variational derivative of the energy functional above.
\begin{align}
  &\operatorname{div}\bm{\sigma} = 0,\quad \bm{\sigma} = \left(g(c)+\eta\right)\frac{\partial W_0}{\partial\bm{\varepsilon}}, \nonumber \\
  & 0 = g'(c) \mathcal{H}^+ - G_c \left[2\xi \Delta c + \frac{w'(c)}{2\xi}\right], \quad \text{where } \mathcal{H}^+:= \max_{\tau\in [0,t]} W_0^+ (\bm{\varepsilon}(x,t))\nonumber\\
  &\forall \bm{x}: W_0^+ < W_0^- \;\Rightarrow c(\bm{x}) = 1 \label{eq:PFFracture}
\end{align}
We replace equation (\ref{eq:PFFracture}b) with a Ginzburg-Landau type evolution law by introducing crack mobility $M$ as
\begin{equation}
  \dot{c} = -M \left[ g'(c) \mathcal{H}^+ - G_c \left(2\xi \Delta c + \frac{w'(c)}{2\xi}\right)\right].\label{eq:PFFractureGL}
\end{equation}
Using equations (\ref{eq:PFFracture}a), (\ref{eq:PFFracture}c), and (\ref{eq:PFFractureGL}), we solve a classic Mode-I fracture propagation problem using a staggered scheme within Alamo. 
The order parameter $\phi_\epsilon$ from the SBM equation (\ref{eq:SBMLinearElasticity}) corresponds to $g(c)$ in the phase fracture equations.
At each time step, as the crack field $c$ evolves, we update the order parameter $\phi_\epsilon$ for the next iteration of the staggered solver. 
\begin{figure}
  \centering
  \begin{subfigure}{0.5\linewidth}
    \includegraphics[width=\linewidth]{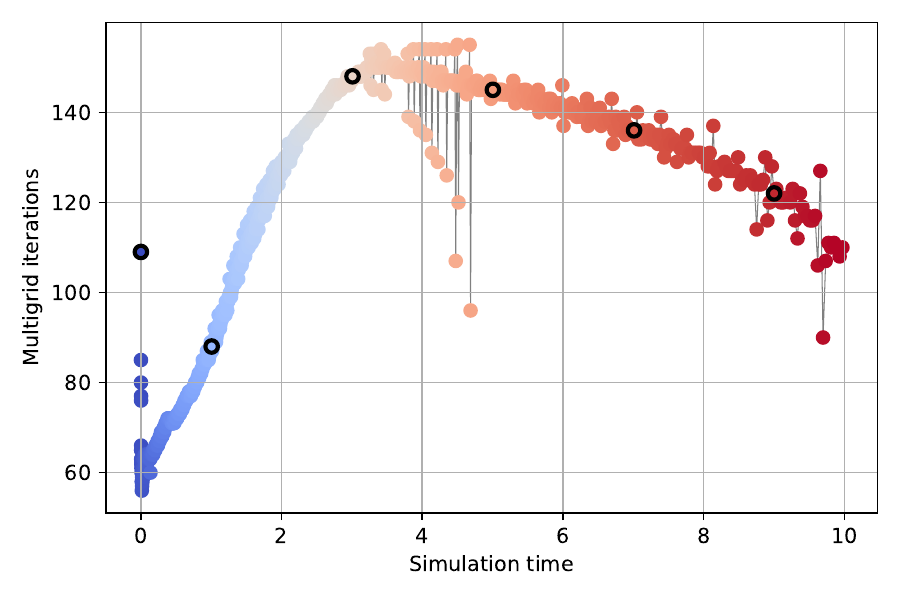}
  \end{subfigure}%
  \begin{subfigure}{0.5\linewidth}
    \includegraphics[width=0.33\linewidth,clip,trim=4cm 0.6cm 0.1cm 0.1cm]{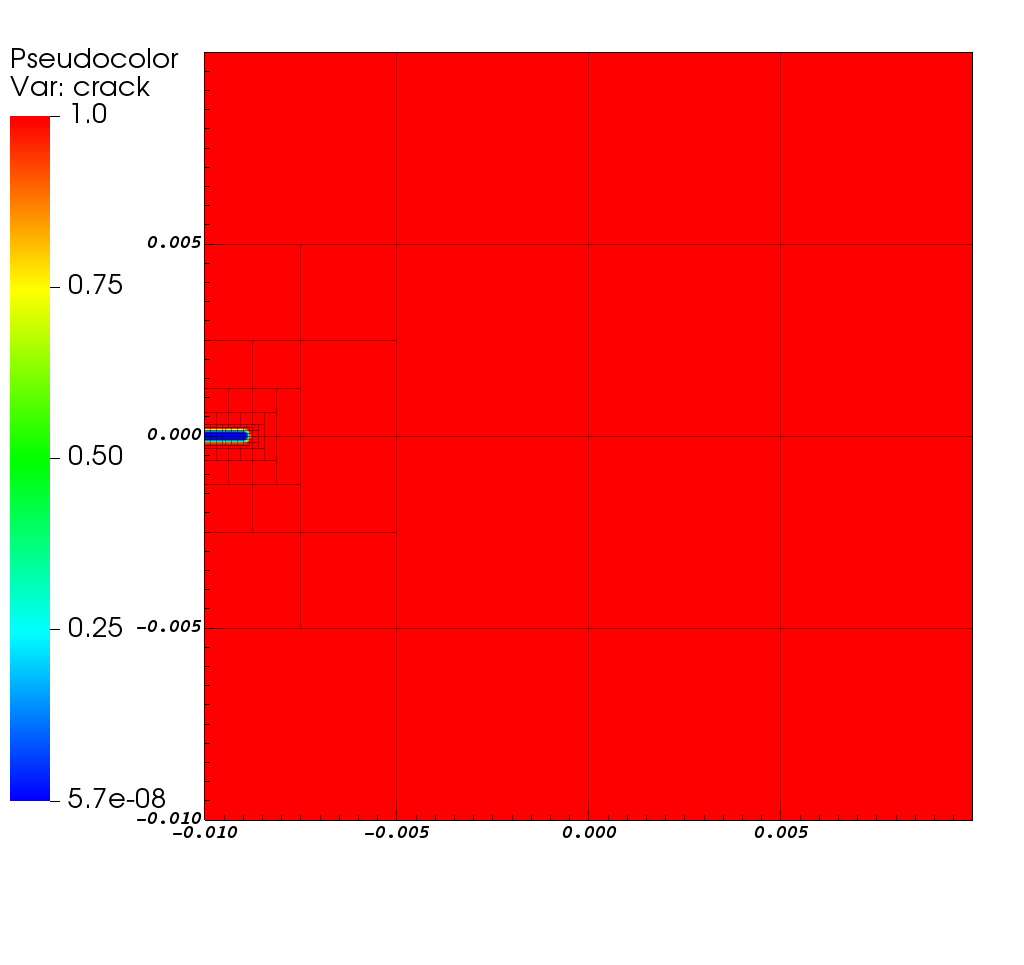}%
    \includegraphics[width=0.33\linewidth,clip,trim=4cm 0.6cm 0.1cm 0.1cm]{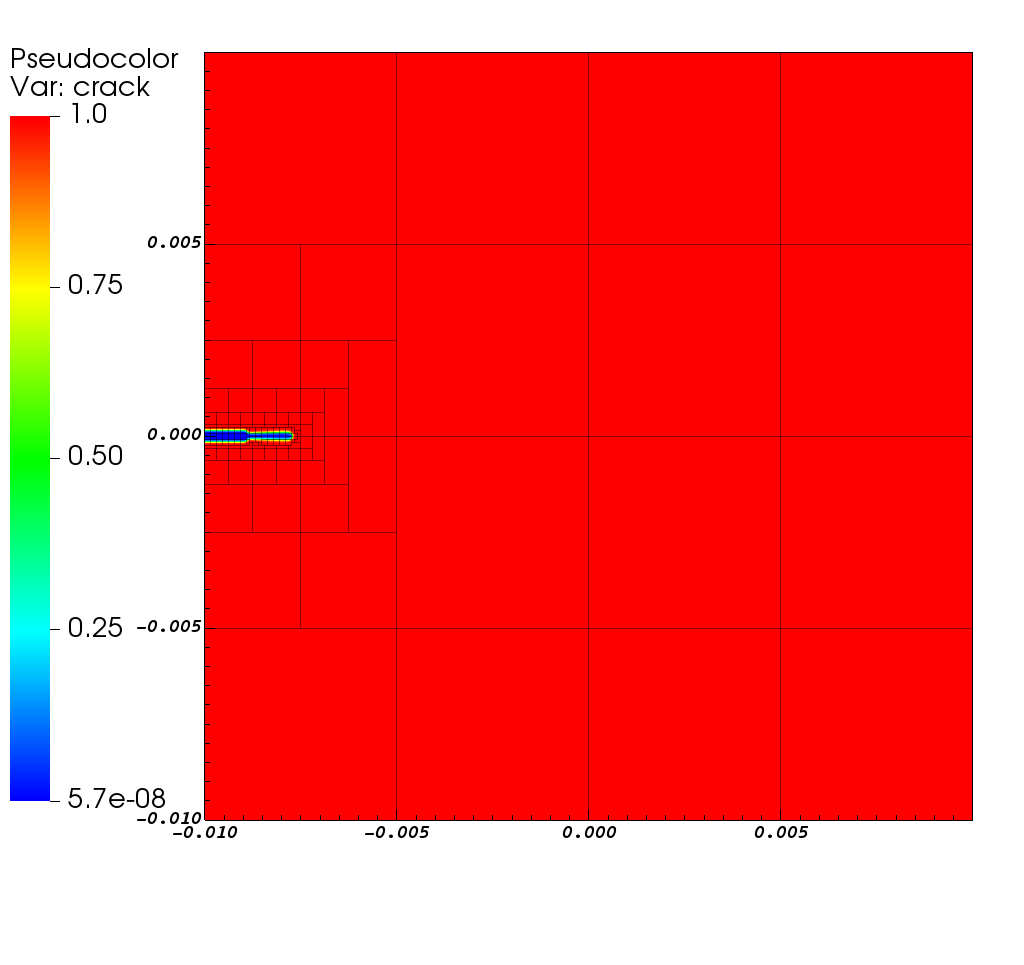}%
    \includegraphics[width=0.33\linewidth,clip,trim=4cm 0.6cm 0.1cm 0.1cm]{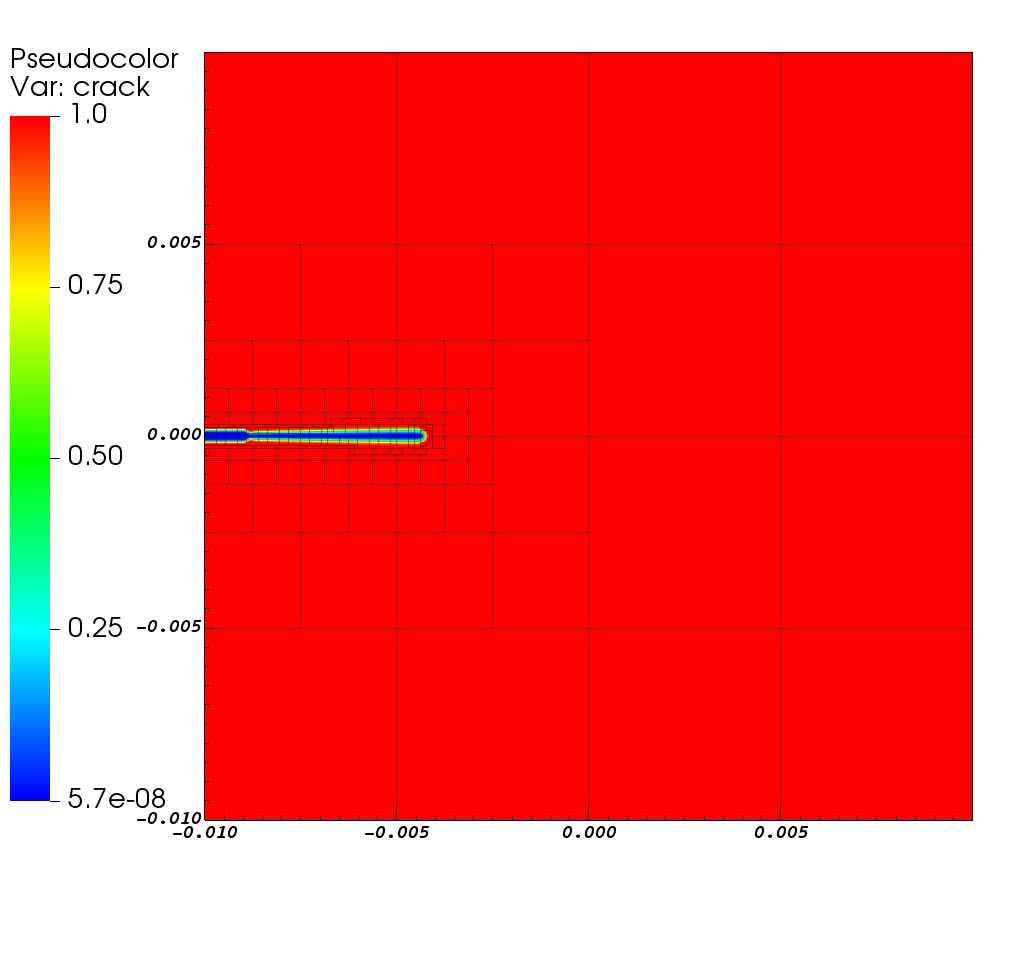}
    \includegraphics[width=0.33\linewidth,clip,trim=4cm 0.6cm 0.1cm 0.1cm]{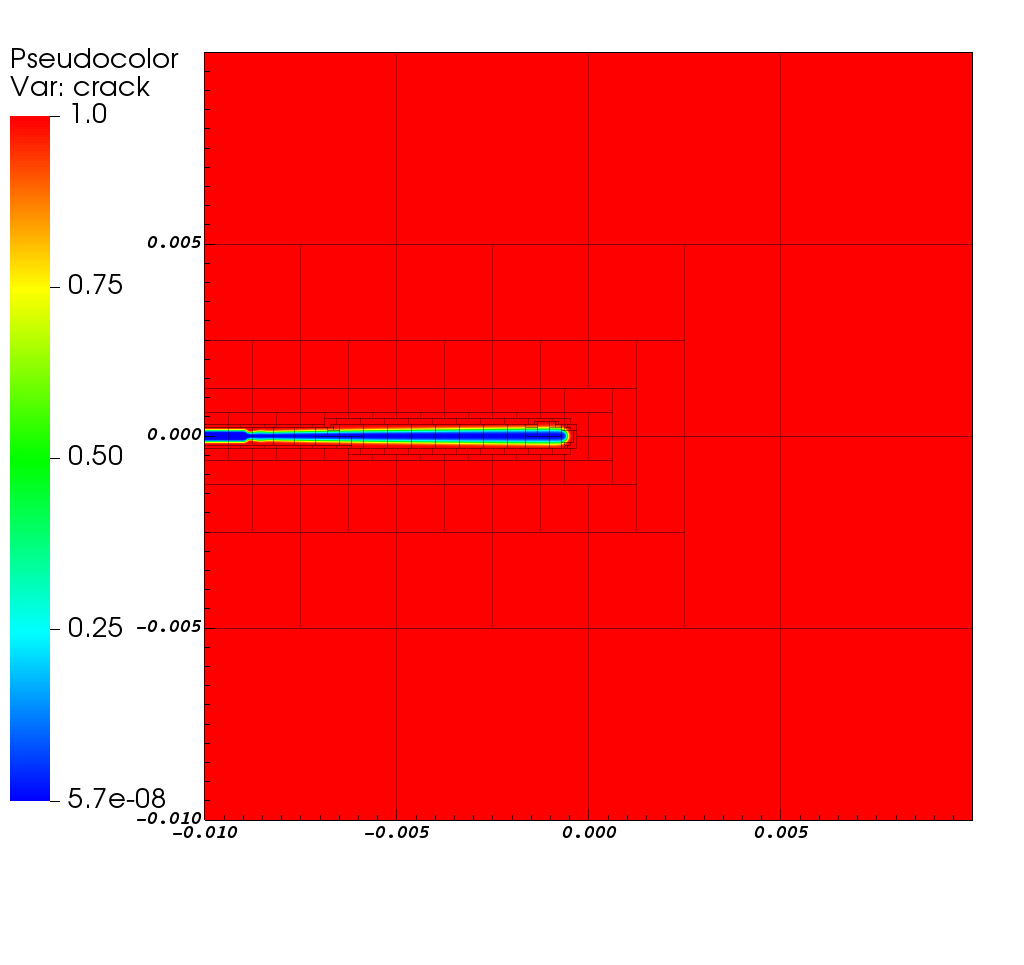}%
    \includegraphics[width=0.33\linewidth,clip,trim=4cm 0.6cm 0.1cm 0.1cm]{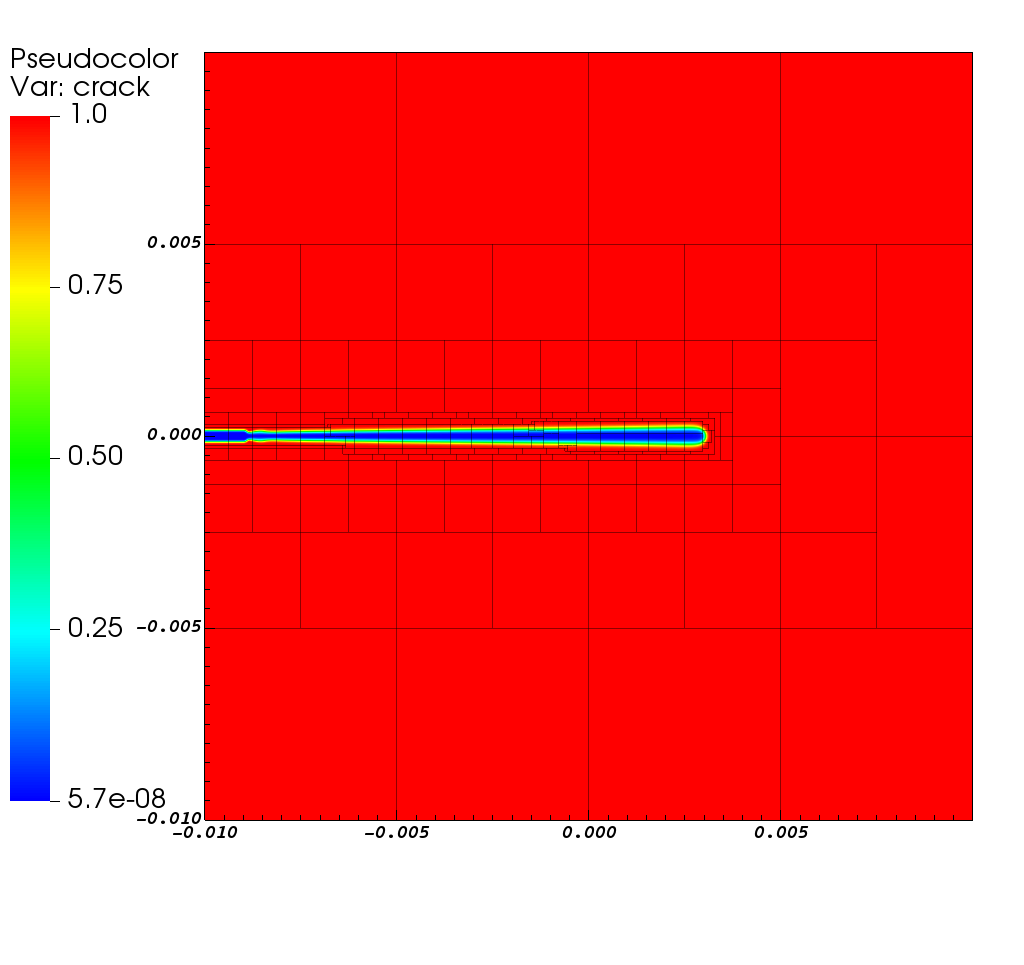}%
    \includegraphics[width=0.33\linewidth,clip,trim=4cm 0.6cm 0.1cm 0.1cm]{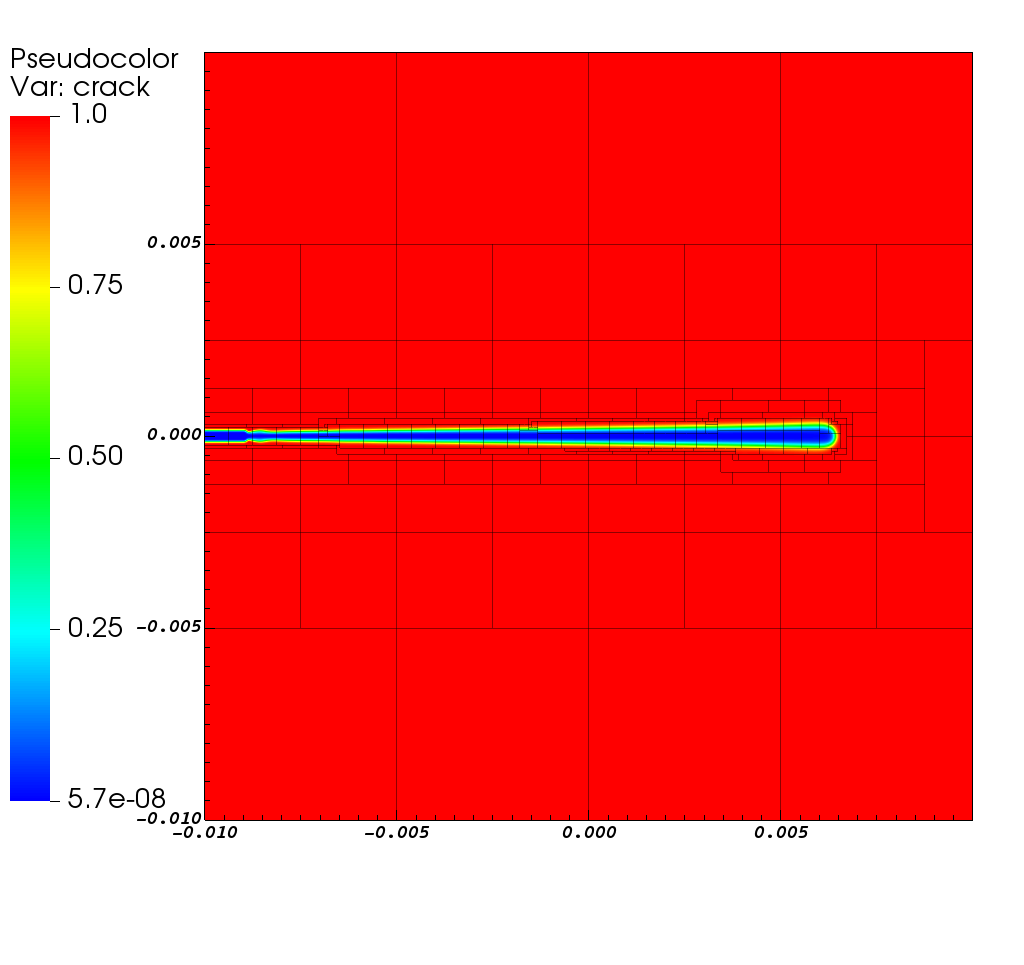}
  \end{subfigure}
  \caption{Phase field fracture - Canonical Mode I loading. (Left) Performance of the MLMG solver during Mode I crack showing the number of MLMG iterations required during simulation time (indicated by position on x axis and by color). (Right) Snapshots showing the phase field evolution at indicated times ($t = 0,1,3,5,7,9$ indicated by black circles) on the performance plot}
  \label{fig:fracture_modeI}
\end{figure}

Figure \ref{fig:fracture_modeI} (right) shows snapshots of mode-I crack propagation over a domain of $\bm{x}\in[-0.01,0.01]\times[-0.01,0.01]$.
We chose a material with $\lambda=121.15 GPa$, $\mu = 80.77GPa$, $\xi=1.0\times 10^{-5}$, $G_c=2700 Pa$, and $M=1.0\times 10^{-5}$.
We initialized a notch of length $1.5\times 10^{-4}$ at the center of the left edge.
We then fixed the bottom boundary, and apply a fixed $y$ displacement of $1.5\times 10^{-5}$ on the top edge. 
We used six levels of refinement on a base grid of $64 \times 64$ to appropriately capture the interface.
The refinement criteria was based on the gradient of the crack field as $|\nabla c|\nabla x| >0.01$.
We performed mesh-regridding every 10 time-steps, with a single time-step being $\Delta t = 10^{-4}$.
As expected, we obtain a steadily propagating crack in the x direction with an adaptively refining grid following the crack field.

We ran this simulation on Auburn University's Easley computing cluster using 32 cores on a single node which took a total of 4.5 hours. 
We note that a major portion of the simulation time was used to perform Ginzberg Landau evolution of the crack field (equation \ref{eq:PFFractureGL}).
Figure \ref{fig:fracture_modeI} (left) shows the number of multigrid iterations needed by the near-singular solver to solve equation (\ref{eq:PFFracture}a).
The solver required 108 iterations for the first elastic solve.
Since the solver uses the previous solution as a starting point, the number of iterations sharply declined immediately after the first elastic solve.
We note a steady increase in required iterations as the crack progresses followed by a peak and slow decline.
We attribute this trend to the changing nature of the mesh as the crack propagates and the fraction of the near-singular domain.
Overall the solver never took more than 160 iterations throughout the entire simulation.

\begin{figure}[t]
  \centering
  \begin{subfigure}{0.5\linewidth}
    \includegraphics[width=\linewidth]{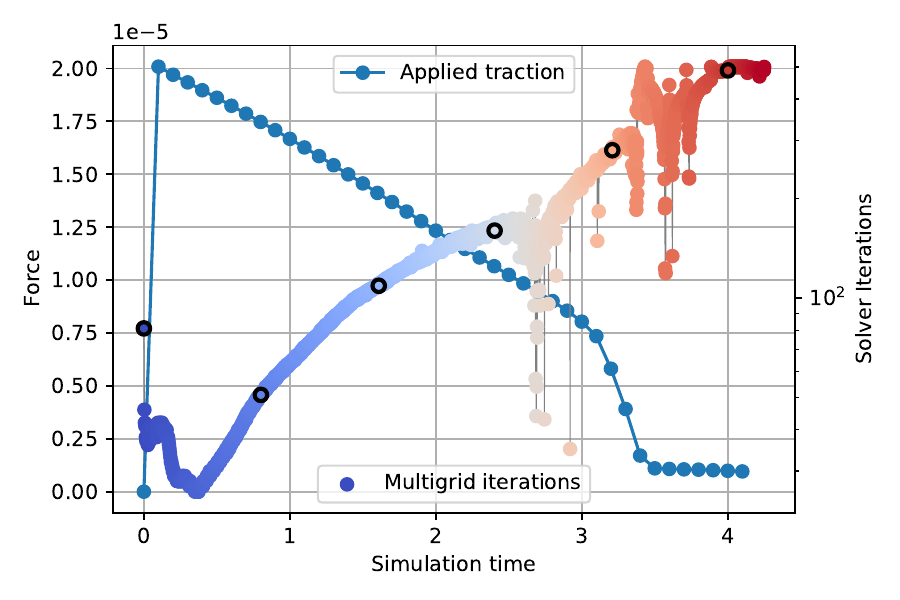}
  \end{subfigure}%
  \begin{subfigure}{0.5\linewidth}
    \includegraphics[width=0.33\linewidth]{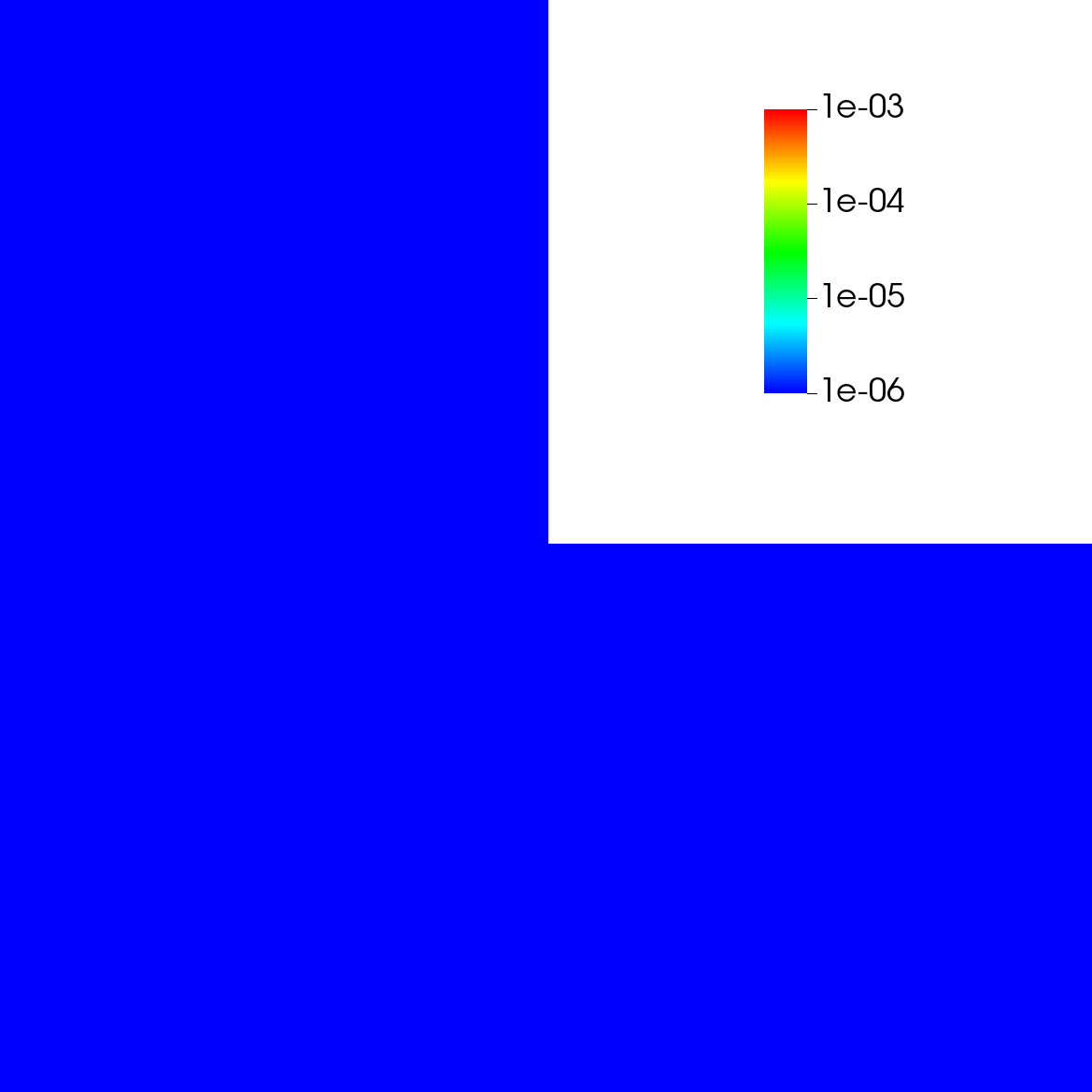}\hspace{1pt}%
    \includegraphics[width=0.33\linewidth]{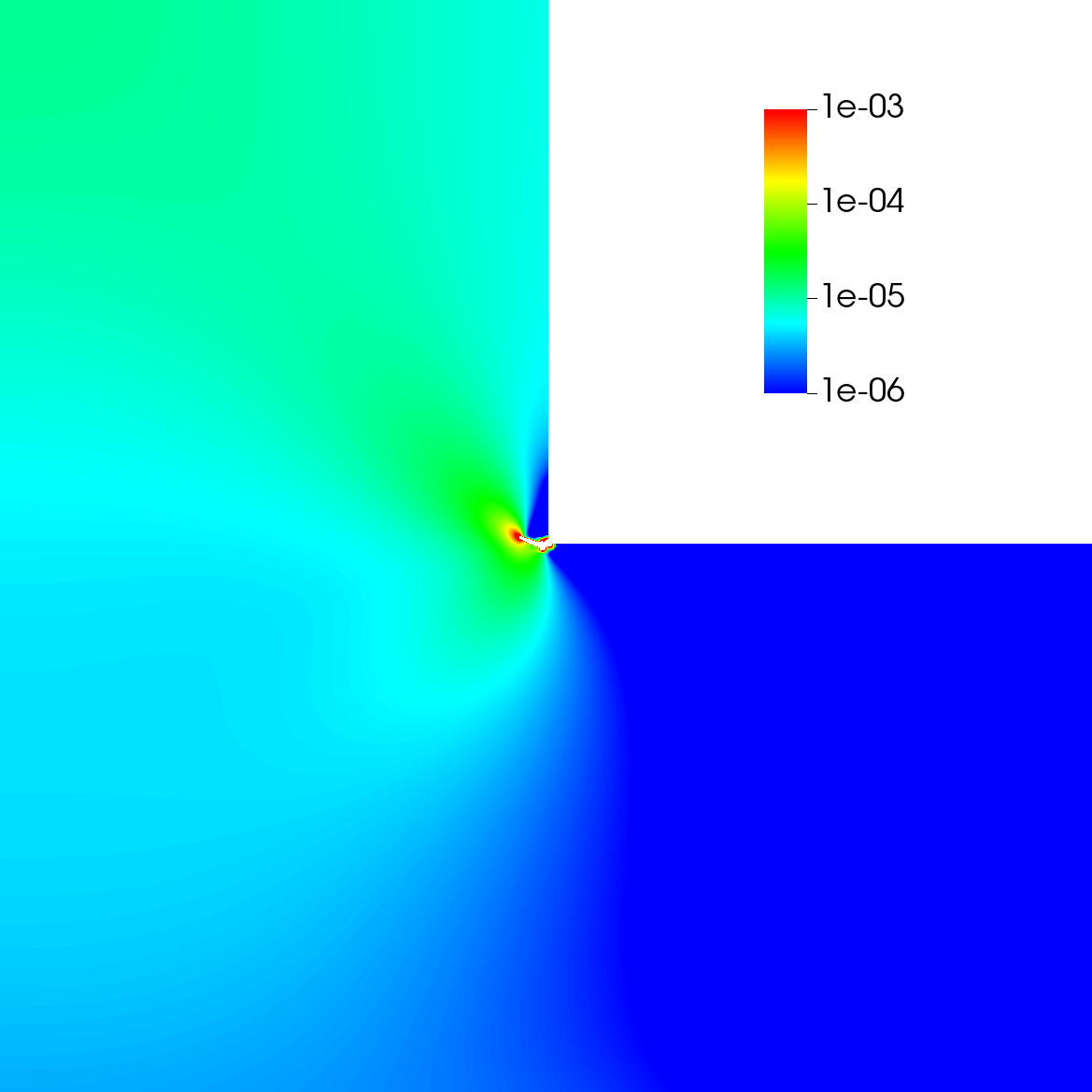}\hspace{1pt}%
    \includegraphics[width=0.33\linewidth]{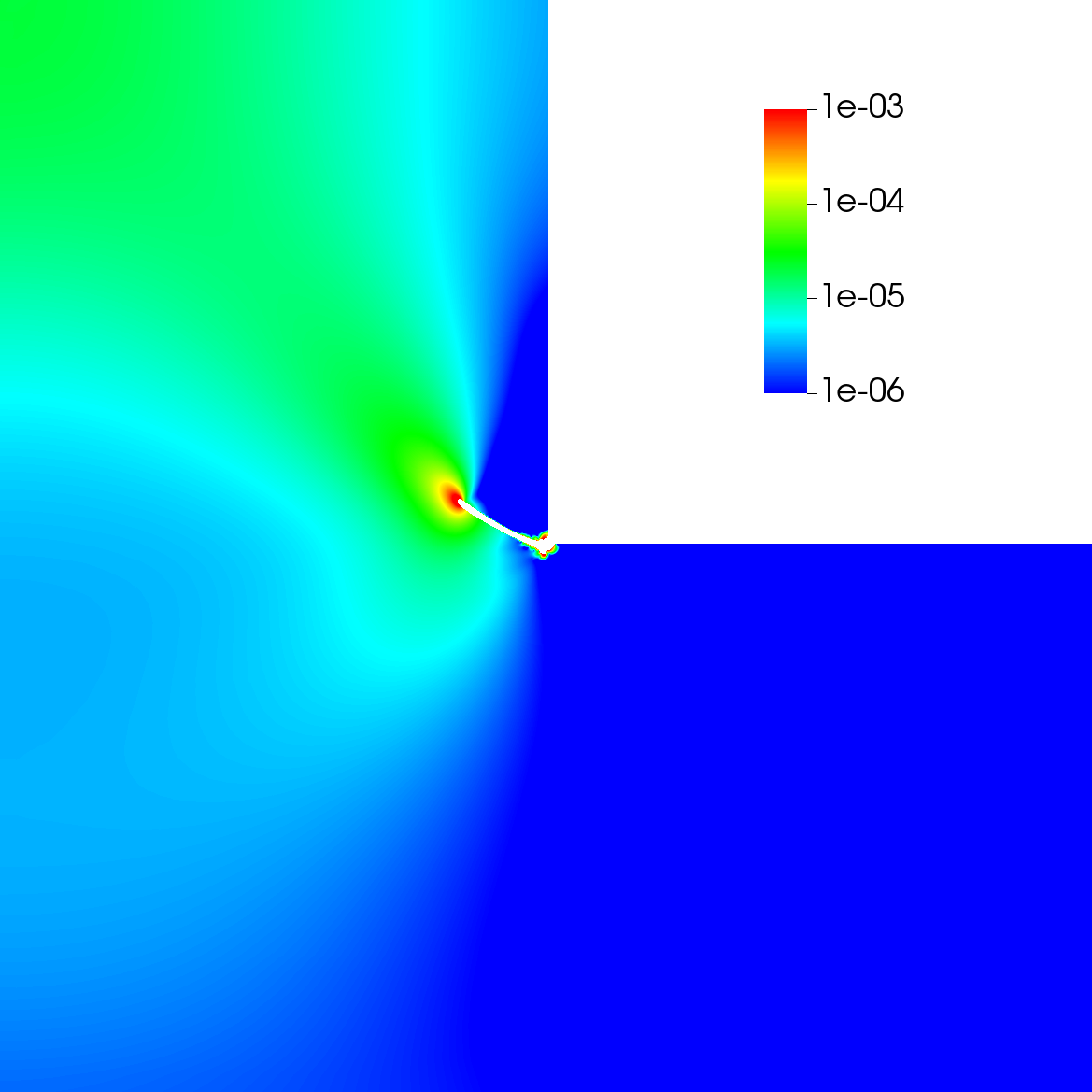}
    \includegraphics[width=0.33\linewidth]{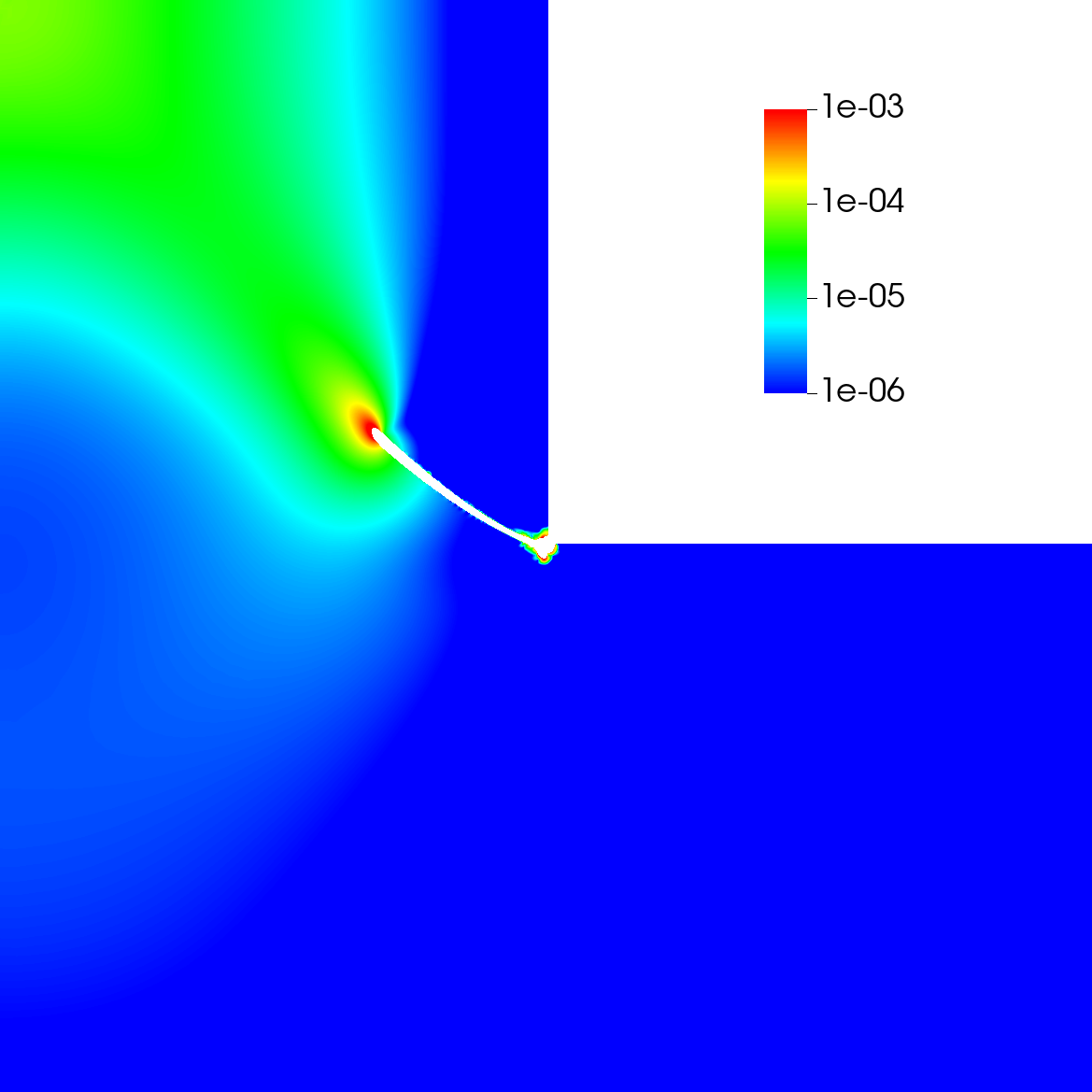}\hspace{1pt}%
    \includegraphics[width=0.33\linewidth]{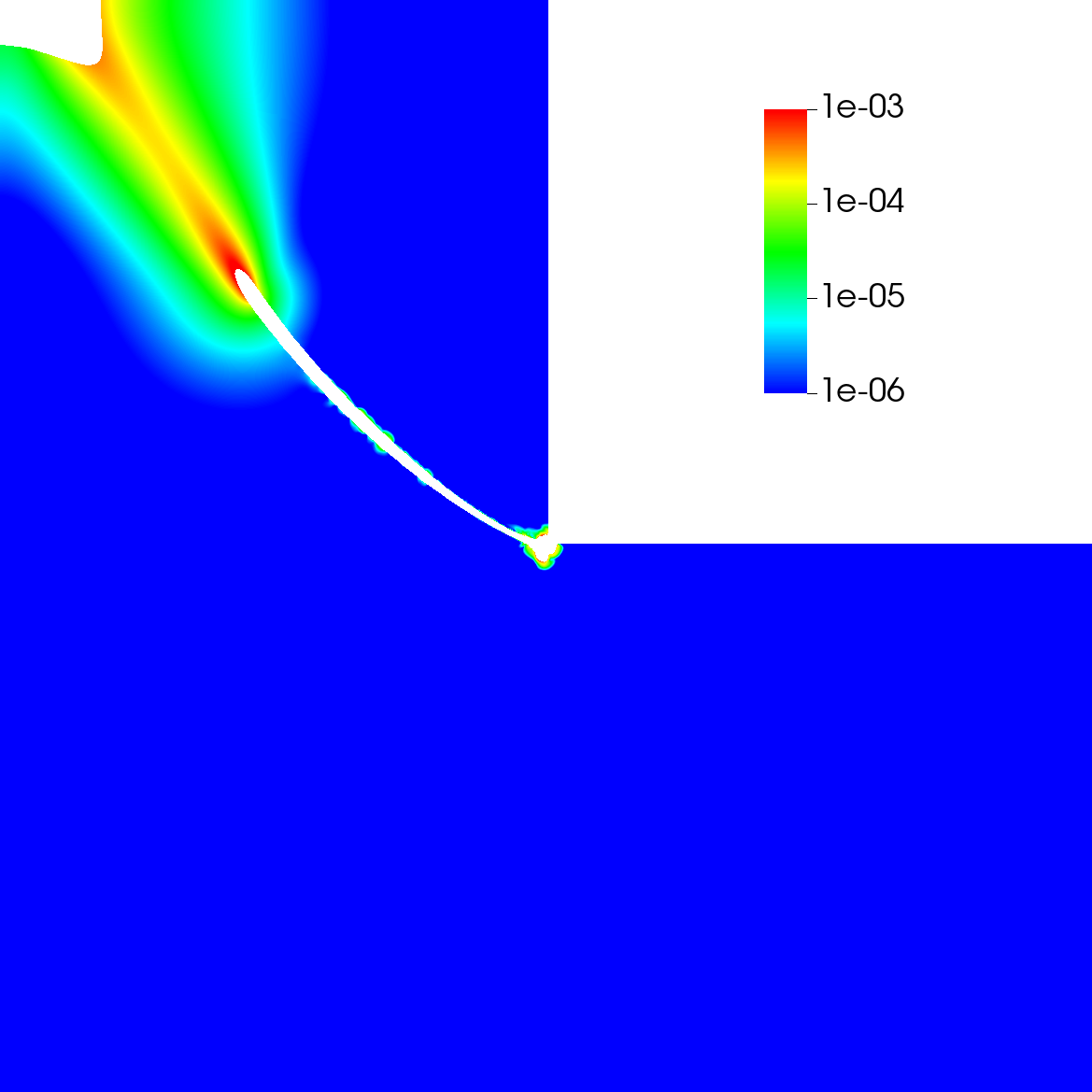}\hspace{1pt}%
    \includegraphics[width=0.33\linewidth]{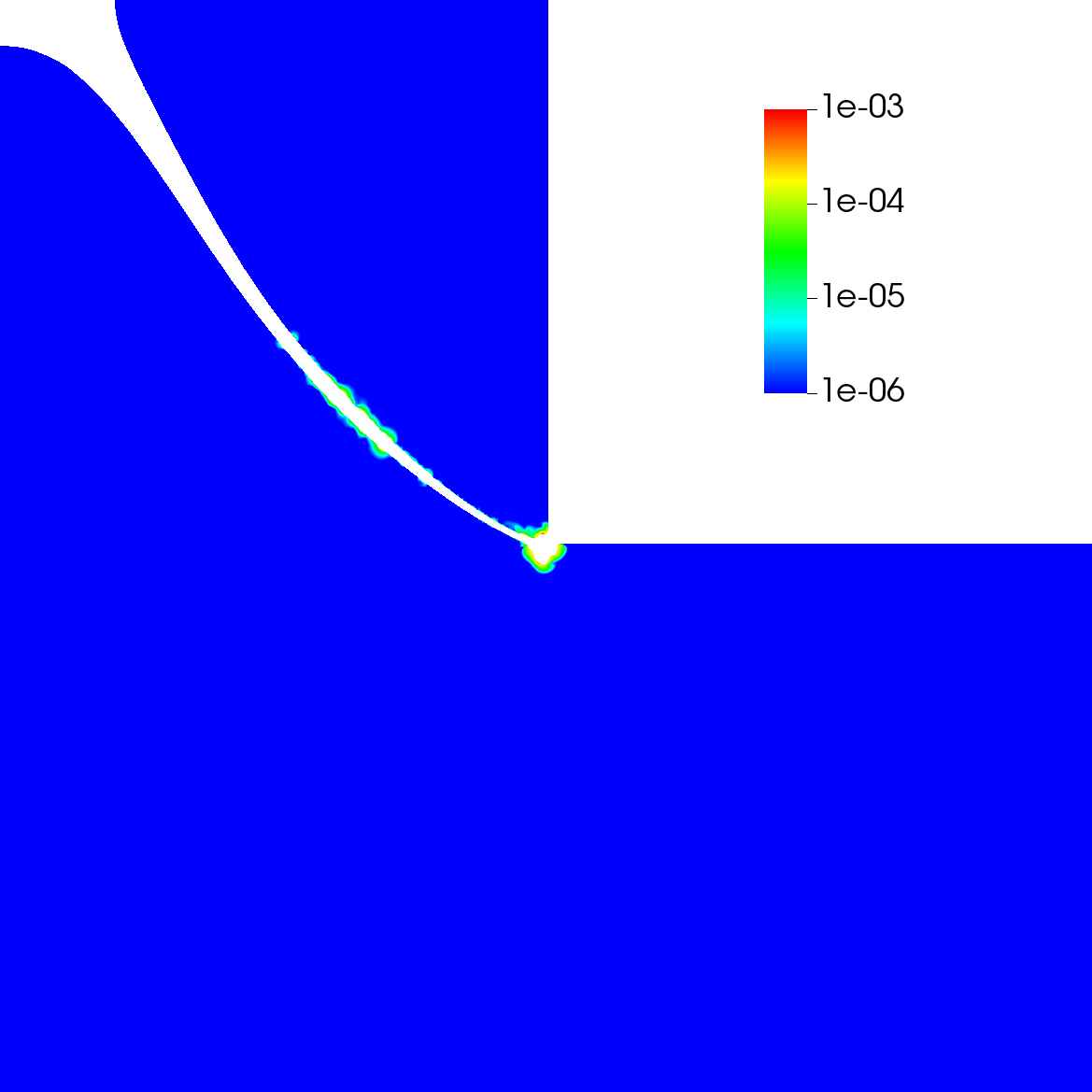}
  \end{subfigure}
  \caption{
    Phase field fracture - Crack nucleated by stress concentration.
    (Left) Performance of the MLMG solver along with the measured applied traction due to the imposed displacement as a function of simulation time (indicated by position on x axis and by color in the iteration plot).
    (Right) Snapshots of the stress state in time ($t=0,0.8,1.6,2.4,3.2,4.0$) as the crack propagates.
    The region where $\psi<0.1$ is colored white.
  }
  \label{fig:crack_lshape}
\end{figure}

We further illustrate the performance of the near-singular solver by studying crack nucleation and propagation due to stress concentration in an L-shaped domain.
We initialized the domain $\Omega := \bm{x}\in [-0.01, 0.01] \times [-0.01, 0.01]$ using a smooth differentiable function with length scale $4\times 10^{-5}$ that takes value $0$ in $\Omega_1 := \bm{x} \in [0, 0.01] \times [0, 0.01]$ and $1$ in $\Omega\setminus \Omega_1$.
We fixed the bottom edge and applied a constant displacement of $1.5\times 10^{-5}$ on the top edge in the y direction.
\Cref{fig:crack_lshape} (right) shows the propagation of crack along with Von-Mises stress distribution in the domain.
As expected, the crack nucleated at the corner $\bm{x}=(0,0)$ with the highest stress concentration and propagated upwards towards the top free surface.
This is confirmed by the \Cref{fig:crack_lshape} (left) where we plot the force (in non-dimensional units) on the top edge. 
We observe a linear decline in the force as the crack propagates, indicating the weakening of the material.
Eventually a secondary crack nucleates at the top left corner, which coalesces with the primary crack causing the final failure of the material.
This is indicated by a sharp decline in measured traction and the snapshot of crack at $t=3.2$.

We performed this simulation on the UCCS INCLINE high-performance computing cluster using 128 cores on a single node. 
The simulation took a total of 6.5 hours most of it, once again, was the Ginzberg Landau evolution of the crack field.
We observe an increase in solver iterations after an initial decline.
The maximum number of iterations required was 500, while the smallest was 48.
Once again, we attribute this pattern to the evolving crack field and near-complete failure of the material.
Overall, we observed results as expected and the solver performed well even near complete failure.

\subsection{Structural topology optimization}
Topology optimization refers to the computational method of determining the geometry of a material or set of materials that produce the optimal result subject to constraints.
Topology optimization generally implies the minimization over a very high dimensional space, the space of all admissible geometries.
Topology optimization has been applied to myriad fields of study, ranging from battery design \cite{mo2021topology} to fluid-structure interaction \cite{andreasen2013topology}.
Structural topology optimization refers specifically to the problem of designing load-bearing structures, subject to constraints typically on the amount of material allowed in a certain volume, that minimizes compliance and maximizes the stiffness of the structure.
Topology optimization has existed as a popular field of study for more than thirty years, stemming from ideas originally proposed more than 150 years ago \cite{logo2020milestones}; today, topology optimization is an entire sub-discipline in its own right.
Topology optimization methods have even found their way into some commercial codes and consequently experienced accelerating usage, partly due to the recent interest in additive manufacturing.

There are a number of prevailing methods for solving structural topology optimization problems.
Common to all structural topology optimization methods is (i) the need to solve the stress equilibrium problem, and (ii) the ability to resolve arbitrary geometry, without {\it a priori} knowledge, with resolution sufficient to resolve lengthscales of interest, and without excessive computational cost.
Following the seminal work by Bendsoe \cite{bendsoe1989optimal,bendsoe1999material}, other methods have included shape derivatives \cite{sokolowski1999topological}, the level set method \cite{allaire2002level,allaire2004structural,wang2003level}, and evolutionary methods \cite{xie1993simple} as techniques for solving the optimization problem.
Recently surging interest in machine learning has led to artificial intelligence-based newcomers, such as generative adversarial networks, that are not necessarily based in physics but are nonetheless capable of generating optimal or near-optimal structures very quickly \cite{li2019non,rade2021algorithmically,chi2021universal}.
A full review of structural topological optimization is well outside this work's scope, so we refer the reader to \cite{sigmund2013topology,jihong2021review} for a more comprehensive overview.

The phase field method is yet another option for solving the structural topology optimization problem.
It is, in some sense, a natural choice, as the phase field method is used specifically for problems involving variable topology.
Phase field was applied to topology optimization by \cite{bourdin2003design}, where the free energy functional contains the elastic strain energy for the given configuration.
Volumetric constraints can be applied by using a conservative Cahn-Hilliard equation that preserves volume \cite{wallin2012optimal}, or by a constrained gradient descent method \cite{wang2004phase,burger2006phase}.
Jeong {\it et al.} implemented volume constraints, as well as additional design constraints, using augmented Lagrange multipliers \cite{jeong2014development}.
A limitation of current phase field methods is the need for high resolution across the diffuse boundary to prevent mesh dependence.
Here, adaptive mesh refinement is needed; while AMR has been applied to phase field topology optimization \cite{salazar2018adaptive,jung2021phase}, the application has been limited.

The finite element method is used nearly universally in topology optimization.
However, the smoothed boundary method for solving near-singular problems presented in this work is ideally suited for solving phase field topology optimization problems.
The method's ability to rapidly regrid, and efficiently solve the near-singular mechanical equilibrium equation, ideally suit it for this application.
In this section, we present results for a basic phase field topology optimization problem.

We use $\eta$ as the order parameter to represent the topology over a domain $\Omega$.
Next, we define the free energy in terms of $\eta$ (where square brackets implicitly indicate a functional over the value and its derivatives) to be:
\begin{align}\label{eq:topop_objective}
  W[\eta] = \int_\Omega\Big(\alpha\eta^2(1-\eta)^2 + \frac{\beta}{2}|\nabla\eta|^2)\Big)\,dx +  \inf_{\bm{u}}\Big[\frac{1}{2}\int_\Omega \nabla\bm{u}\cdot ((\eta+\zeta)^2\mathbb{C})\,\nabla\bm{u}\,dx - \int_{\partial_2\Omega}\bm{u}\cdot\bm{t}^0\,dx \Big],
\end{align}
subject to the constraint
\begin{align}\label{eq:topop_volume_constraint}
  \int_\Omega \,\eta\,dx = V_0
\end{align}
where $\alpha$ and $\beta$ are numerical parameters controlling segregation and boundary energy, $\mathbb{C}$ is the fourth order elasticity tensor, $\bm{t}_0$ is prescribed surface traction, and $V_0$ is the allowable volume of material.
(We note that $\alpha\sim\frac{1}{\epsilon}, \beta\sim\epsilon$ where $\epsilon$ controls the boundary width. 
We retained $\alpha$ and $\beta$ for simplicity.)
As discussed above, the volume requirement induces a constraint on the optimization problem.
We adopt a straightforward regularization, allowing the optimum to be found by a modified version of the Allen-Cahn equation:
\begin{align}\label{eq:topop_evolution}
  \frac{\partial\eta}{\partial t} = -L\Big(\frac{\delta}{\delta\eta}W + \lambda(t)\Big(\int_\Omega\,\eta\,dx - V_0\Big)\Big)
\end{align}
where $L$ is a mobility parameter and $\partial/\partial\eta$ is the variational derivative.
The function $\lambda(t)$ is a Lagrange multiplier that enforces the constraint in \cref{eq:topop_volume_constraint}.
In general, we let $\lambda(t)$ tend towards infinity as $t\to\infty$, to allow the constraint to be approached at a gradual rate.
In problems with multiple optima, the form for $\lambda(t)$ may determine which local minimum is found.
Each evaluation of \cref{eq:topop_evolution} requires the evaluation of the elastic minimization problem in \cref{eq:topop_objective}, which is solved using the proposed method.
As in all examples, $\eta$ is a cell-centered field while the displacement is node-centered.
Node-to-cell averages are computed for each iteration.

\begin{figure}
  \begin{minipage}{0.2\linewidth}
    \includegraphics[height=1.9cm,clip,trim=0cm 0 22cm 0]{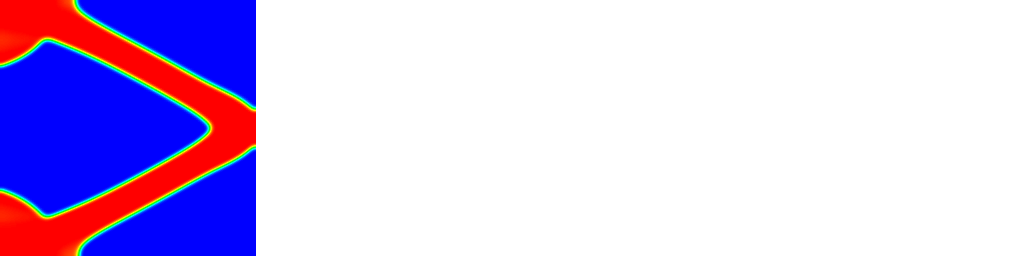}
    \includegraphics[height=1.9cm,clip,trim=0cm 0 22cm 0]{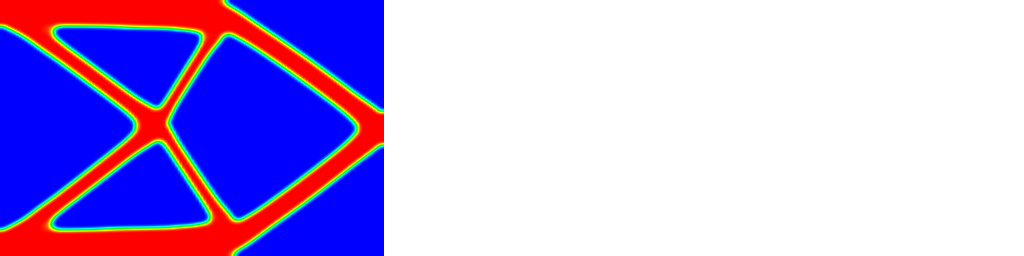}
  \end{minipage}%
  \begin{minipage}{0.33\linewidth}
    \includegraphics[height=1.9cm,clip,trim=0cm 0 0cm 0]{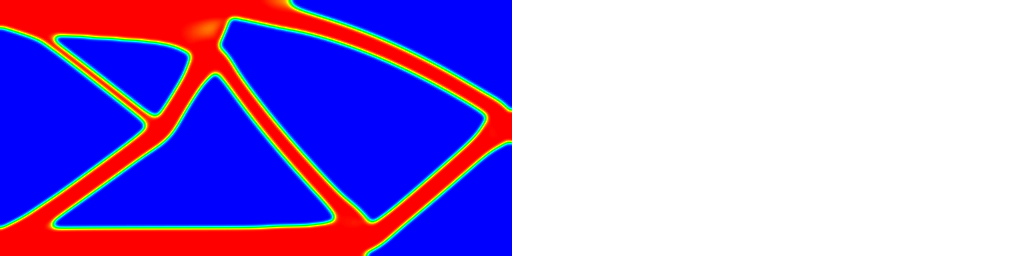}
    \includegraphics[height=1.9cm]{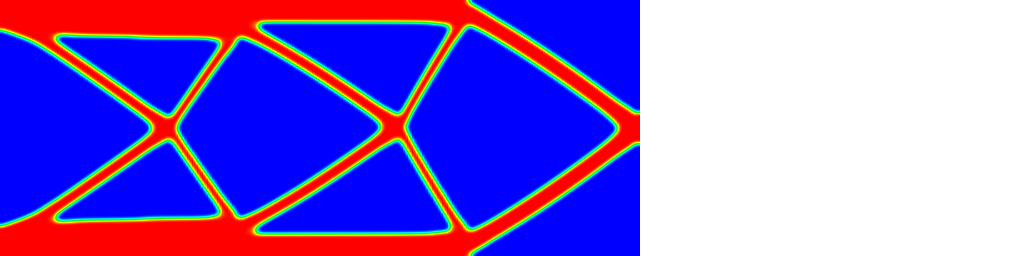}
  \end{minipage}%
  \begin{minipage}{0.33\linewidth}
    \includegraphics[height=1.9cm]{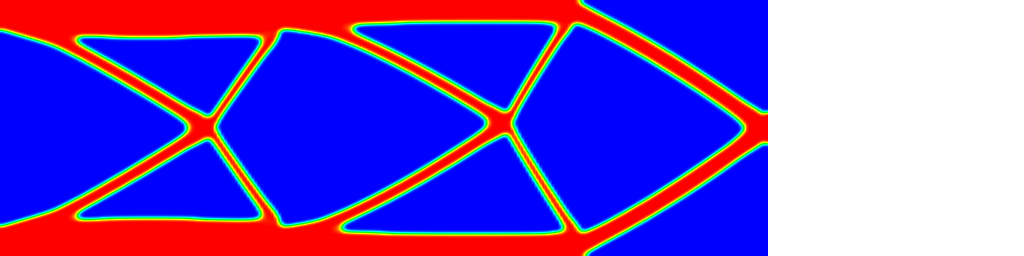}
    \includegraphics[height=1.9cm]{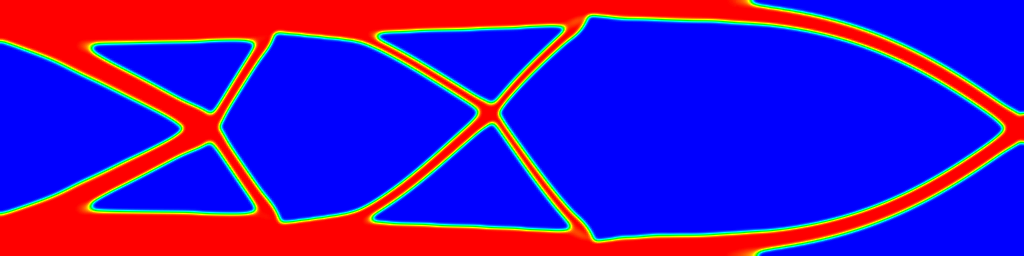}
  \end{minipage}
  \caption{
    Topology optimization results for point load supported by a cantilever, wtih 25\% fill and increasing aspect ratio.
  }
  \label{fig:topop_thumbnails}
\end{figure}

We considered the classical problem of a load supported by a cantilevered structure.
We selected dimensionless values for all quantities.
For the two-dimensional results, we chose the domain to be $1$ in height, and ranging from $1$ to $4$ in length.
We used base-level grids of 32x32, 64x32, and 128x32 corresponding to the different aspect ratios.
We used the powers of two for all grid dimensions, at the cost of some non-square unit cells, to optimize the multigrid solver's performance.
We used a total of three AMR levels, with a refinement criterion $|\nabla\eta||\Delta\bm{x}| \ge 0.05$.
We used a less restrictive refinement criterion for the strain, $\varepsilon$, although it did not contribute significantly to the results.
We performed mesh regridding at each phase field iteration step.

For the cases considered here, we constrained the volume to 25\% of the domain volume.
We chose the phase field parameters as $\alpha=200,\beta=0.01$, $\lambda=400.0$.
We used an isotropic material model with $E=1480$, and $\nu=0.22$.
We applied a point load of magnitude $-0.1$ per unit length at the center of the right face over a region $0.01$ in height.
We specified Neumann conditions for the order parameter, and it is possible to see slight artifacts at the edges of the domain that result from this condition.
We chose the regularization parameter $\zeta$ to be $0.01$.
While we considered smaller values, they did not affect the performance of the solver.
However, they did affect the nucleation behavior of the phase field method, causing spurious material segments to be generated and eventually resulting in some instability that caused the solution to diverge.
Therefore, we attribute this to limitations of the phase field model and leave further optimization of the model to future work.

The algorithm's results were generally as expected for models of this type (\Cref{fig:topop_thumbnails}).
For the 1x1 case, the result is a fairly simple triangular brace structure.
Varying the volume fraction and boundary width terms generally did not produce a substantial change.
Increasing the domain to 1.5x1 produced a truss-like structure that is generally in line with the canonical result for this standard problem.
Increasing again to 2x1 produced an irregular, asymmetric structure, likely due to the inaccessibility of an optimal symmetric structure for that aspect ratio.
Asymmetries were common in this work, especially for cases with higher boundary penalization terms.
Since we did not initialize the problem asymmetric perturbations, any deviation generally stemmed from perturbations induced by the regridding algorithm.
For the 2.5x1 case, the resulting structure is similar to the 1.5x1 case except with a secondary truss structure in addition to the first.
The top and bottom support beams also decrease in width so as to prioritize support near the wall where the bending moment is the highest.
As the domain continues to increase in length, the double support structure is elongated but does not generate the third structure.
This is due to the width of the diffuse interface, as the thickness of the third support structure would be thinner than the diffuseness of the boundary.
Therefore, it is not possible to nucleate the third structure.
We emphasize that we did not account for buckling in this model. 
The slenderness of the beams would, in reality, result in buckling that renders the structure unstable.
While it is possible to modify the model to avoid this effect by differentiating between tension and compression, it is outside the scope of the present work.

\begin{figure}
  \begin{subfigure}{0.4\linewidth}
    \includegraphics[width=\linewidth]{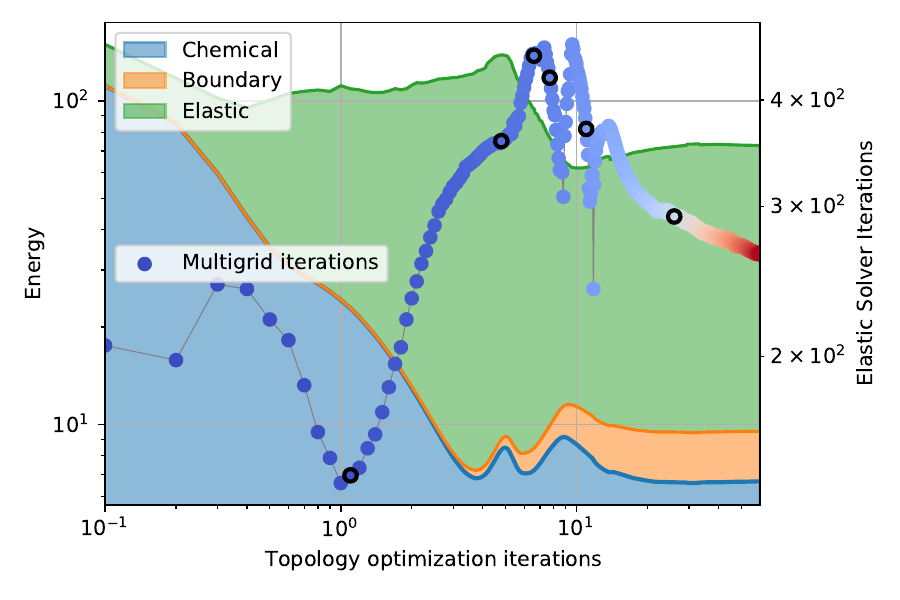}
    \caption{}
    \label{fig:topop_history}
  \end{subfigure}
  \begin{subfigure}{0.6\linewidth}
    \begin{minipage}{0.33\linewidth}\centering
      \includegraphics[width=\linewidth,clip,trim=10cm 0 0 0]{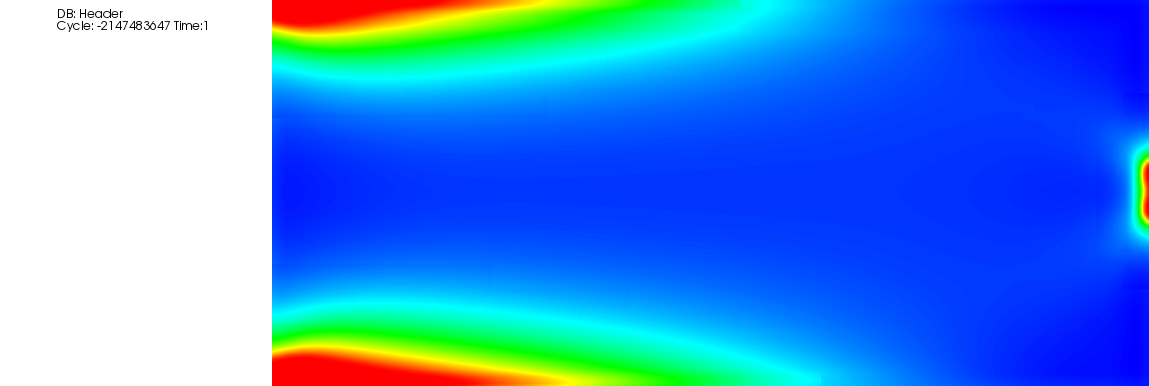}
      $t=1.0$
    \end{minipage}%
    \begin{minipage}{0.33\linewidth}\centering
      \includegraphics[width=\linewidth,clip,trim=10cm 0 0 0]{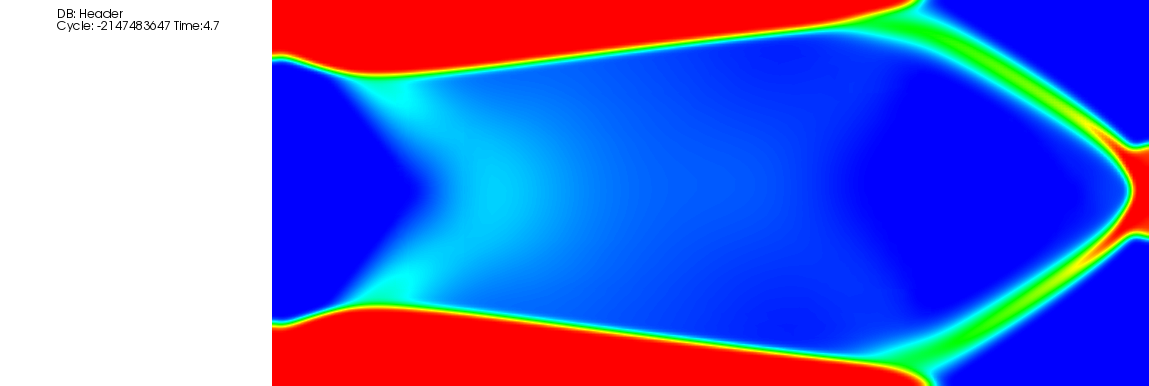}
      $t=4.7$
    \end{minipage}%
    \begin{minipage}{0.33\linewidth}\centering
      \includegraphics[width=\linewidth,clip,trim=10cm 0 0 0]{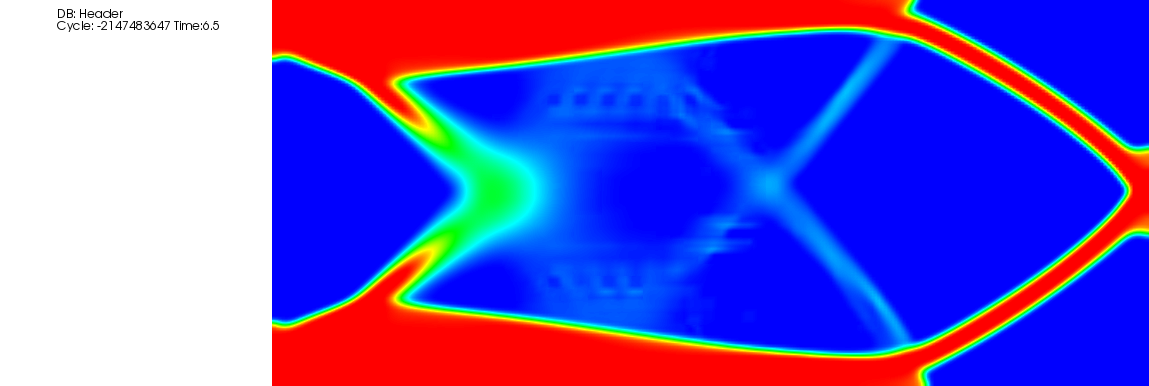}
      $t=6.5$
    \end{minipage}\vspace{10pt}
    \begin{minipage}{0.33\linewidth}\centering
      \includegraphics[width=\linewidth,clip,trim=10cm 0 0 0]{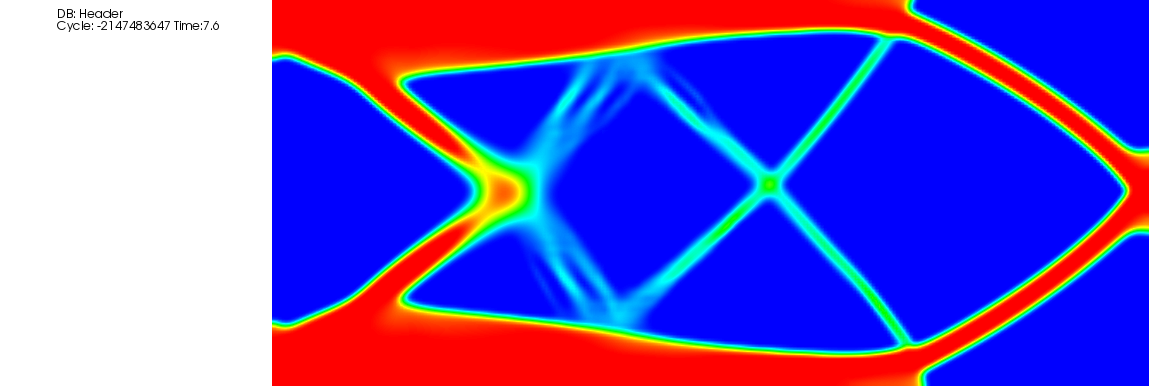}
      $t=7.6$
    \end{minipage}%
    \begin{minipage}{0.33\linewidth}\centering
      \includegraphics[width=\linewidth,clip,trim=10cm 0 0 0]{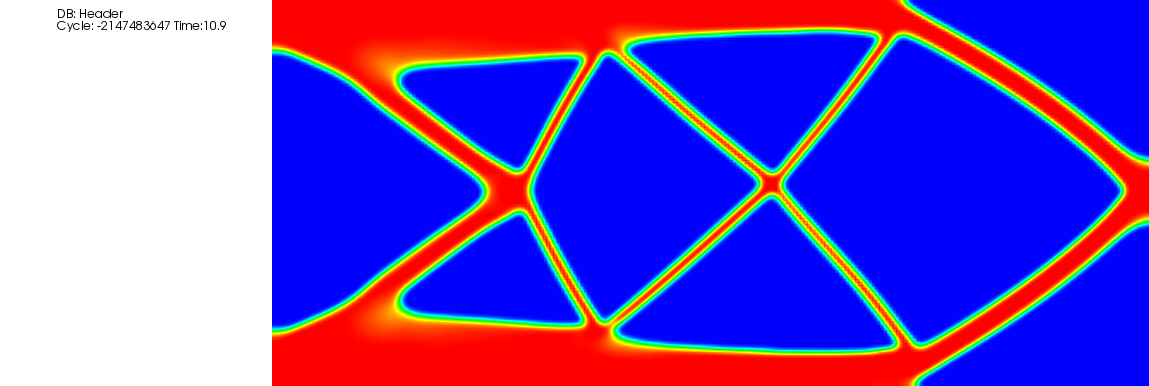}
      $t=10.9$
    \end{minipage}%
    \begin{minipage}{0.33\linewidth}\centering
      \includegraphics[width=\linewidth,clip,trim=10cm 0 0 0]{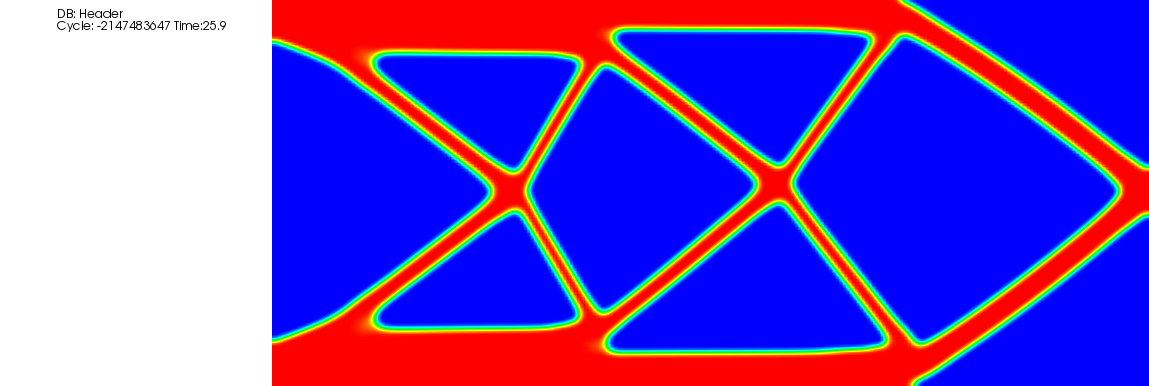}
      $t=25.9$
    \end{minipage}
    \caption{}
    \label{fig:topop_snapshots}
  \end{subfigure}
  \caption{Progression of the energy and the corresponding performance of the elastic solver during topology optimization. (a) Contributions of chemical potential, boundary energy, and elastic energy to the objective function. Number of solver iterations required are indicated by the position along the x axis and by color.  (b) Snapshots corresponding to the black dots in (a).}
  \label{fig:topop_designevolution}
\end{figure}

Of particular interest is the performance of the solver during the solution of the phase field topology optimization problem.
We present the history of the design evolution for a representative structure in \Cref{fig:topop_history}.
The colored regions represent contributions of the energy functional corresponding to the chemical potential, boundary energy, and elastic energy.
(We note that the regions are stacked and are plotted on a log scale in the x and y axes.)
The initial iterations are dominated by the chemical potential, which drives the segregation and, as a result, a sharp increase in the elastic energy.
Eventually, support structures are spontaneously nucleated, rapidly relieving the elastic stress and moving the solution toward its equilibrium state.
We superimpose the performance of the multigrid solver (gray points) and plot it with respect to the right axes, also on a log scale.
In general, the solver always converged linearly within a couple hundred iterations and generally no more than 400.
(See \ref{sec:convergence_data} for more details on solver convergence.)
We notice a decrease in the solver's performance when topological changes occurred, which can be connected to the presence of problematic ``islands'' and ``peninsulas'' in the solution. 
This can be seen by observing the design evolution (\Cref{fig:topop_snapshots}) at the indicated points (black circles) on the performance curve.

All of the simulations took less than an hour to complete on a single node with 32 cores.
We used the UCCS INCLINE cluster to perform the simulations, but only due to the large number of simulations that were considered.
We note that all the 2D results can be quickly reproduced on a desktop computer in a matter of minutes to hours (depending on domain size and computing power).
We emphasize that the objective of this work is not to compete with commercially available topology optimization codes, but to demonstrate the verstility of the method and its ability to easily adapt to a diverse range of problems..

\begin{figure}
  \includegraphics[width=0.33\linewidth]{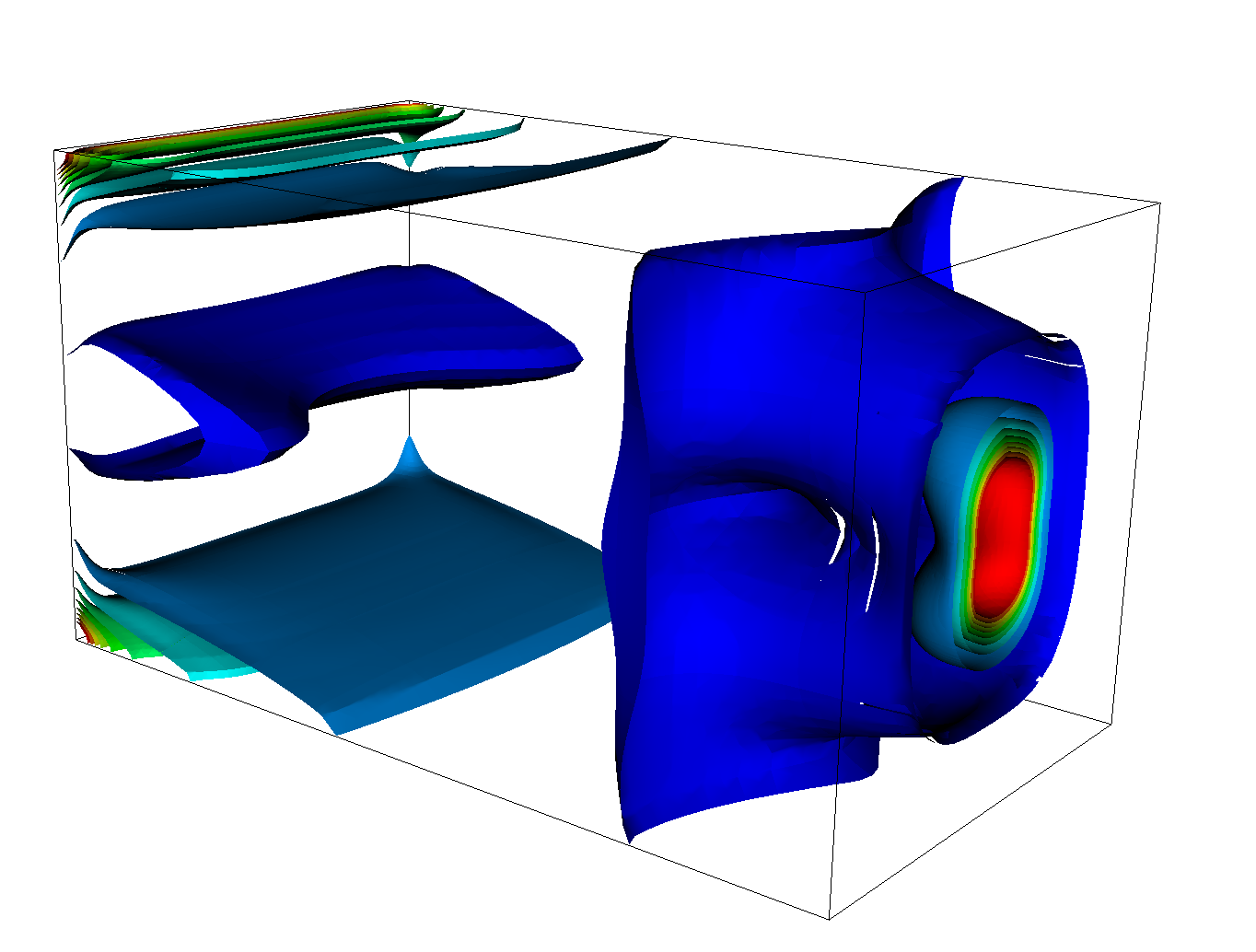}
  \includegraphics[width=0.33\linewidth]{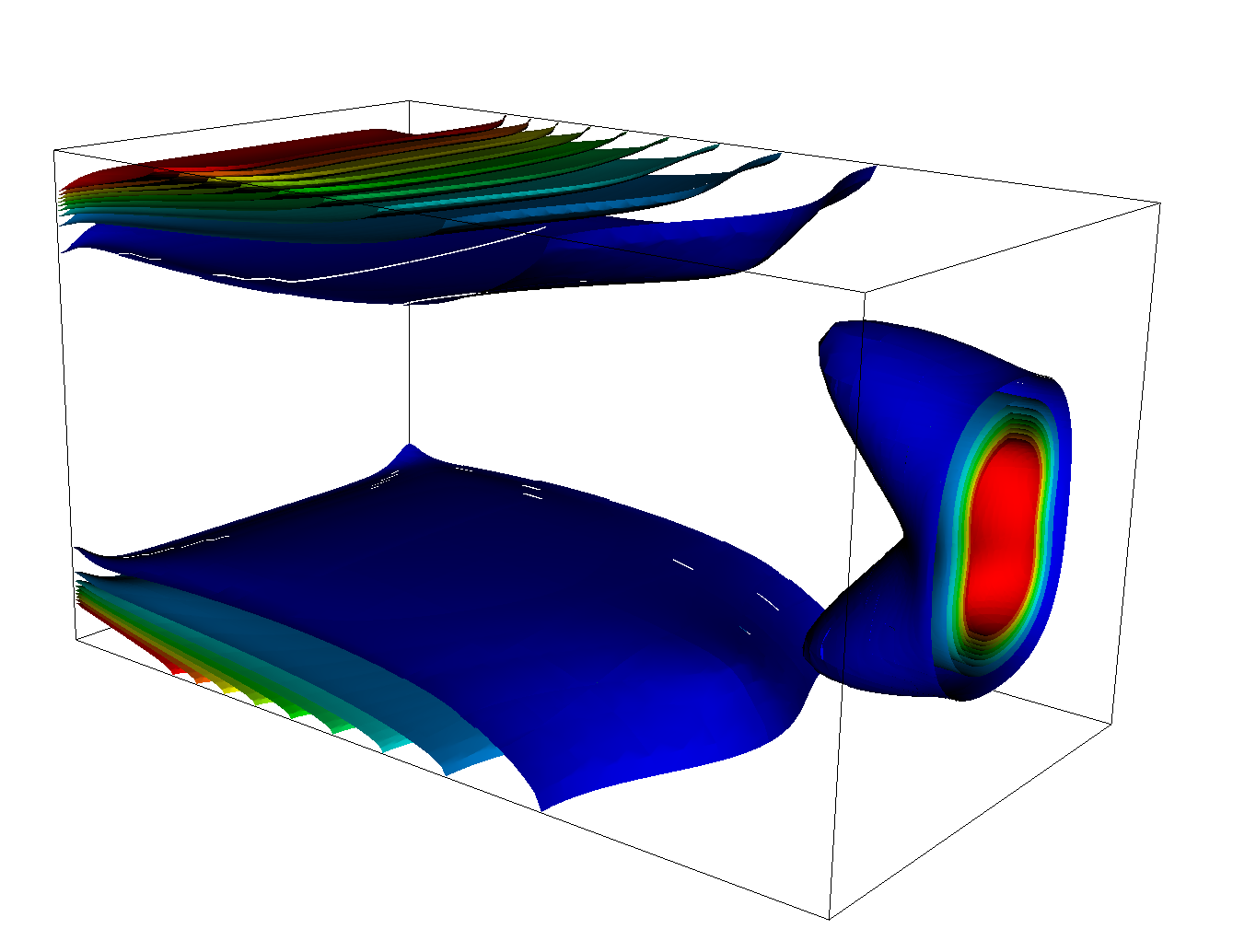}
  \includegraphics[width=0.33\linewidth]{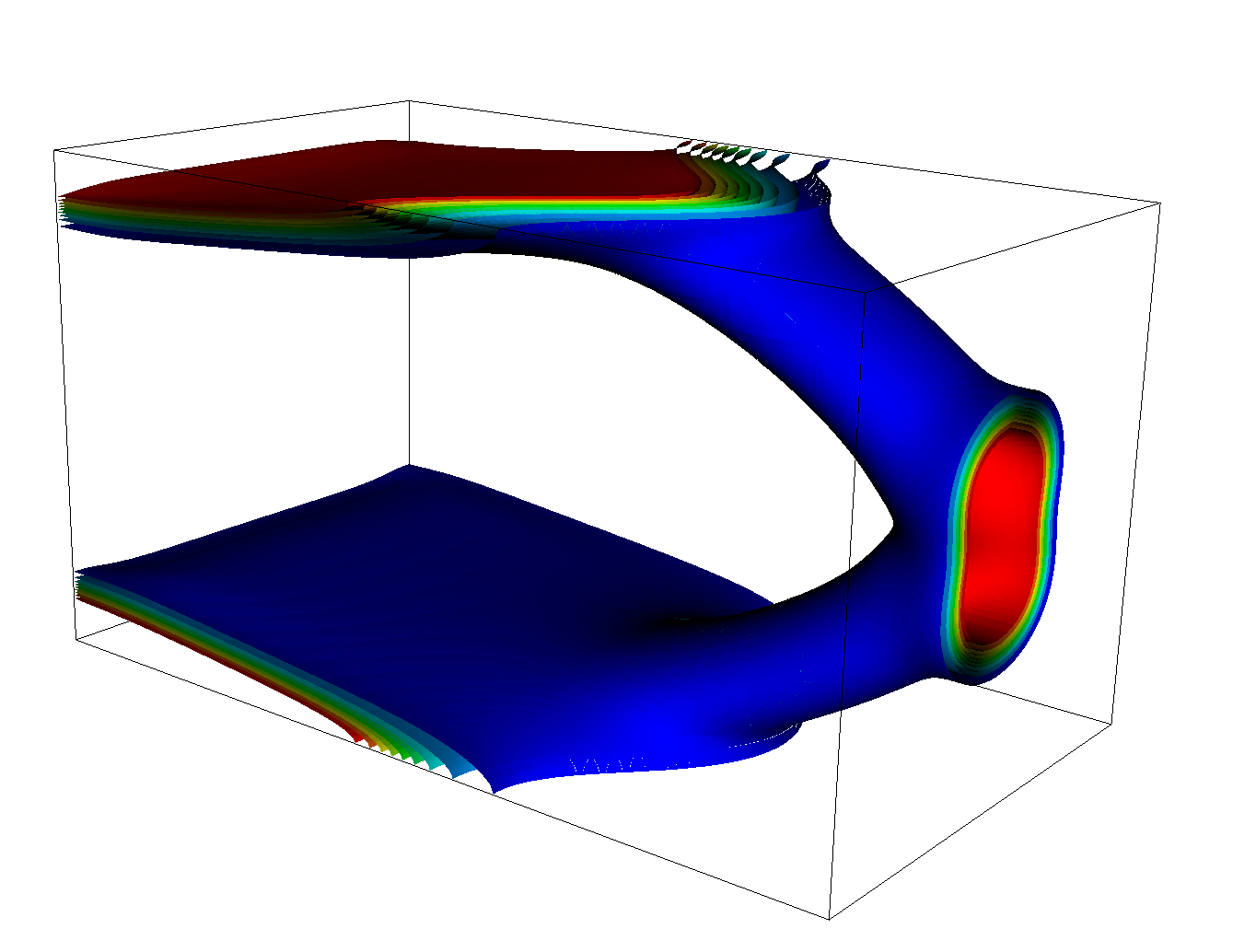}
  \includegraphics[width=0.33\linewidth]{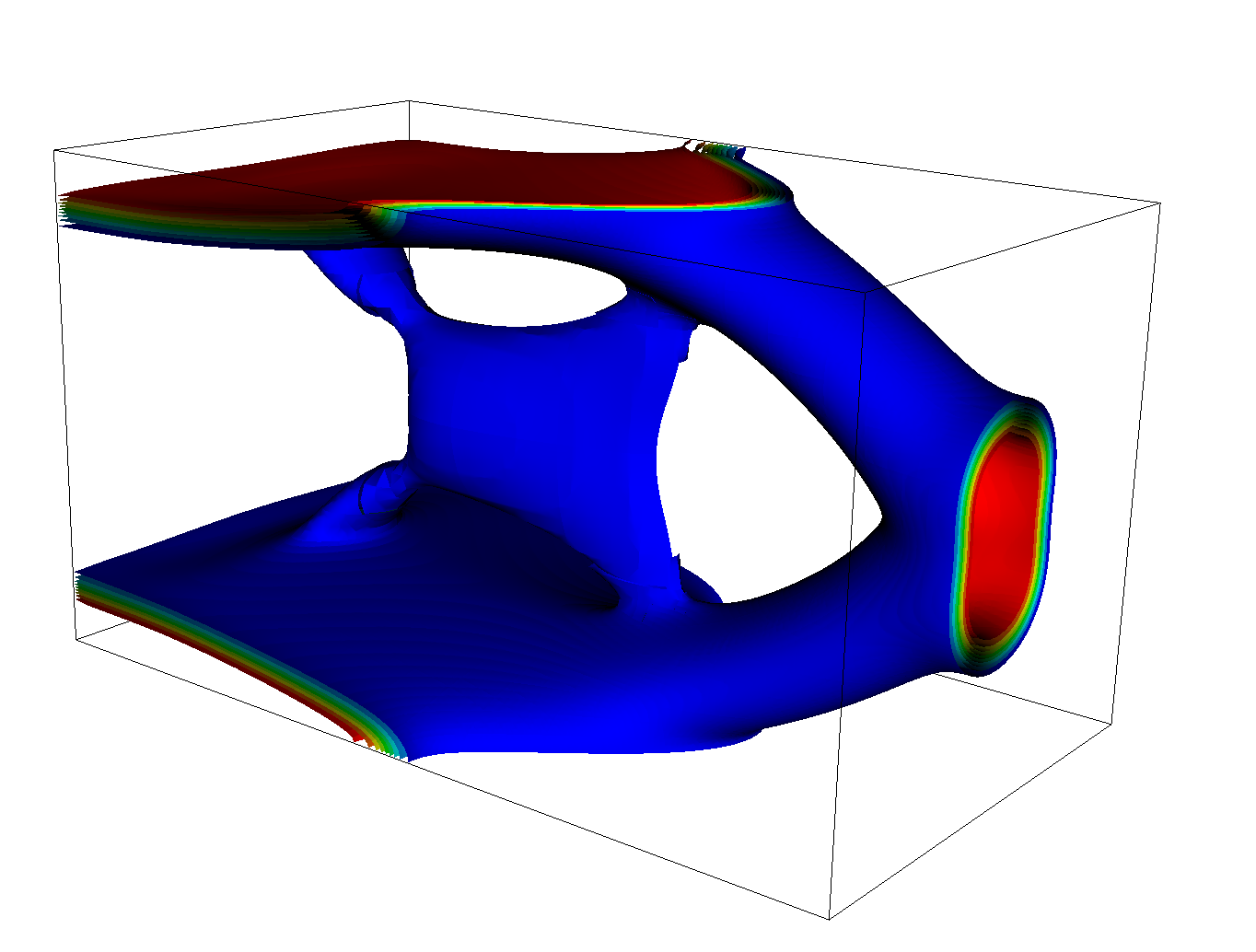}
  \includegraphics[width=0.33\linewidth]{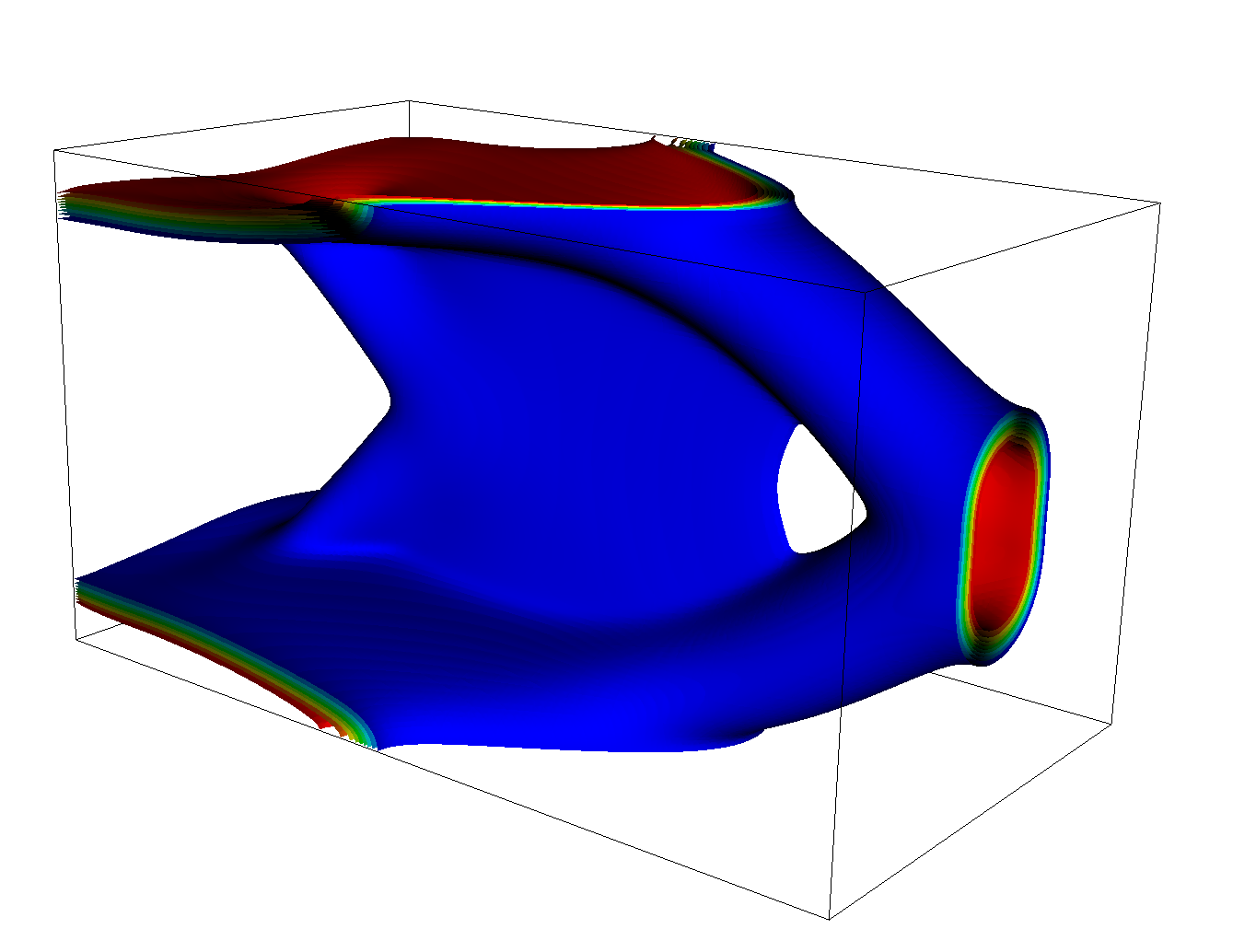}
  \includegraphics[width=0.33\linewidth]{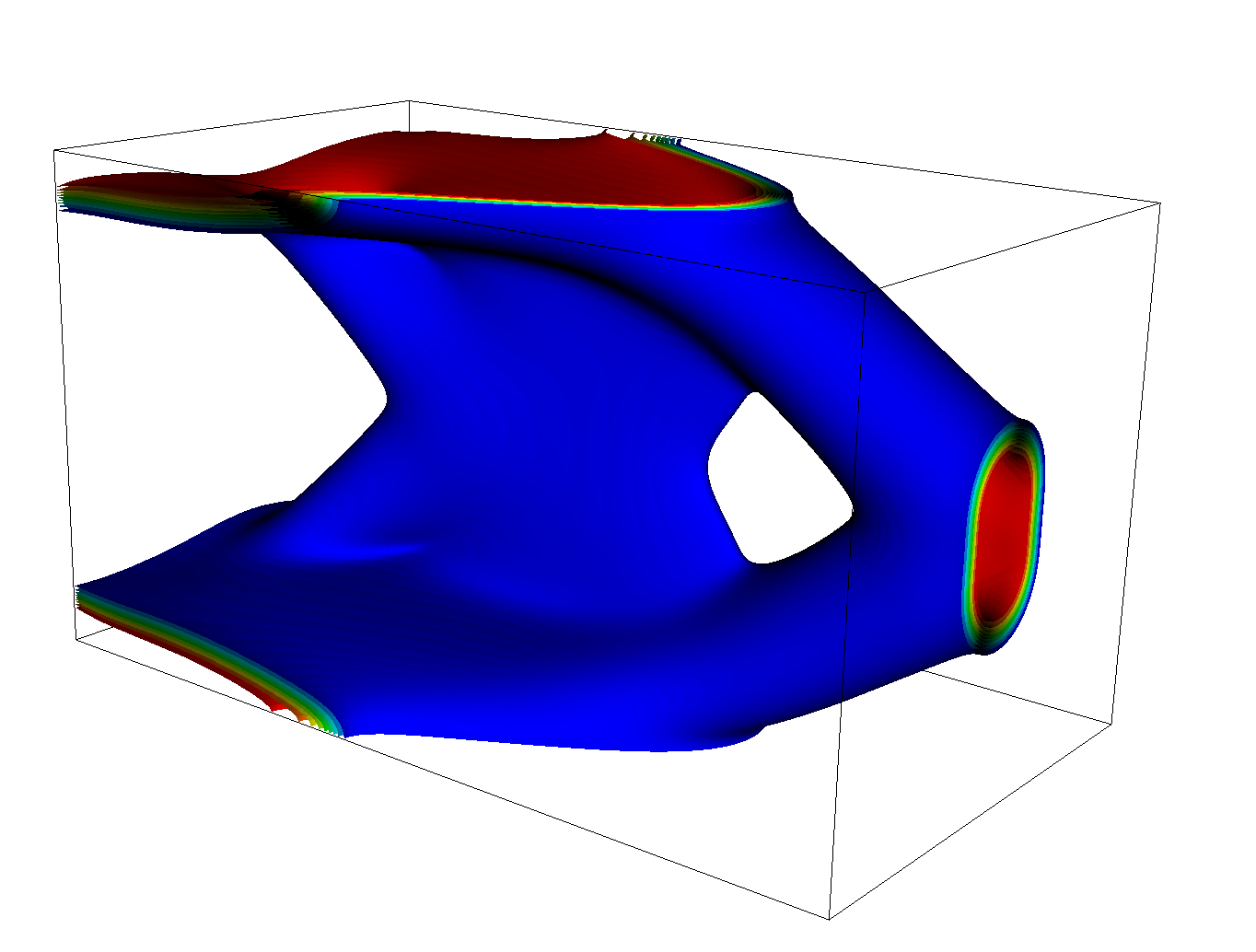}
  \caption{
    Three dimensional topology optimization results for a cantilever load.
    Contours correspond to increments of $\Delta\eta=0.1$, snapshots are shown every 10 timesteps.
    The full simulation took 4 hours to run on 128 processors; the approximate result was reached after about 1 hour of computing time.
    Two levels of refinement are used for a total of 3 BSAMR levels.
  }
  \label{fig:topop_3d}
\end{figure}

Finally, we tested the method in three dimensions (Figure \ref{fig:topop_3d}).
We chose a configuration that was generally similar to the 2D case except for the following differences.
We chose the domain as 1.8x1.0x1.0, with a base grid of 64x32x32.
We used a larger value of $\beta=0.1$.
We applied the same point load at the right end, except that it was applied at the center over an area of $0.1\times 0.1$.
The figure shows the evolution of $\eta$ during the solution, with isosurfaces plotted at increments of $\Delta\eta=0.1$.
We ran the simulation on a single node (128 processors) of the INCLINE cluster for four hours, although we observed that the result converged well before the conclusion of the run time (about an hour).
As with the two-dimensional case, the result is a truss-like structure with a central support mechanism.
The result differs from the 2D case in that the central truss is replaced with a webbed structure.
(Once again, we note that our model did not account for buckling, which would significantly change the final result.)
In general, we observed results as expected and satisfactory performance from the model in 3D. 

\section{Conclusion}
In this work, we presented a comprehensive computational approach for solving near-singular problems with the smoothed boundary method and block-structured adaptive mesh refinement.
Problems in solid mechanics are nearly universally solved using weak form methods (finite elements), and this work aims to formalize and elucidate the theory and best practices associated with solving problems of interest (specifically, near-singular problems) using the near-singular, smoothed boundary, strong-form approach.
In particular, a key insight of this work is the essential role that is played by placing certain quantities at node or cell centers in finite difference elasticity.
This has been well-known in fluid mechanics and has generally been enforced implicitly in solid mechanics due to the prevalence of strong-form methods.
This, combined with the model vector space construction and the reflux-free method for elasticity on BSAMR grids, results in a powerful technique for quickly solving many problems of interest in solid mechanics.
We demonstrated this approach's versatility and the solver's performance by studying three problems: plasticity, fracture, and topology optimization.

We conclude by discussing the limitations of this work.
As with all diffuse boundary methods, a higher resolution is needed than most discrete boundary approaches.
BSAMR is able to alleviate this computational cost significantly, but the number of points will still scale with the amount of surface area, even for surfaces with low stress.
(Discrete boundary methods, by contrast, can afford low resolution at uninteresting boundaries.)
While this may be unavoidable, we believe that the proposed approach offers increased versatility and ease of implementation.
We also note that while the solver performs well on near-singular problems, it can still struggle on problems that exhibit extreme irregularity or high stress concentrations.
In this regard, it is comparable to equivalent solvers for alternative approaches.
Finally, we emphasize that the examples considered here do not reflect the state-of-the-art in the fields of plasticity, fracture, and topology optimization; rather, we have used tested and well-understood models in order to showcase the performance of the elastic solver.
Implementing more sophisticated methods using this solver shall be left to future work.

\section{Acknowledgements}

VA acknowledges the Auburn University Easley Cluster for support of this work.
BR acknowledges support from Lawrence Berkeley National Laboratory, subcontract \#7645776, and from the Office of Naval Research, grant \#N00014-21-1-2113. 
This work used the INCLINE cluster at the University of Colorado Colorado Springs. 
INCLINE is supported by the National Science Foundation, grant \#2017917.
The authors wish to thank Dr. Scott Runnels at the University of Colorado Boulder, whose insights on cell and node-based fields led to key breakthroughs in solver development.

\appendix

\section{Convergence data}\label{sec:convergence_data}

\setcounter{figure}{0}

\begin{figure}[h]
  \centering
  \begin{subfigure}{0.5\linewidth}
    \includegraphics[height=5.5cm]{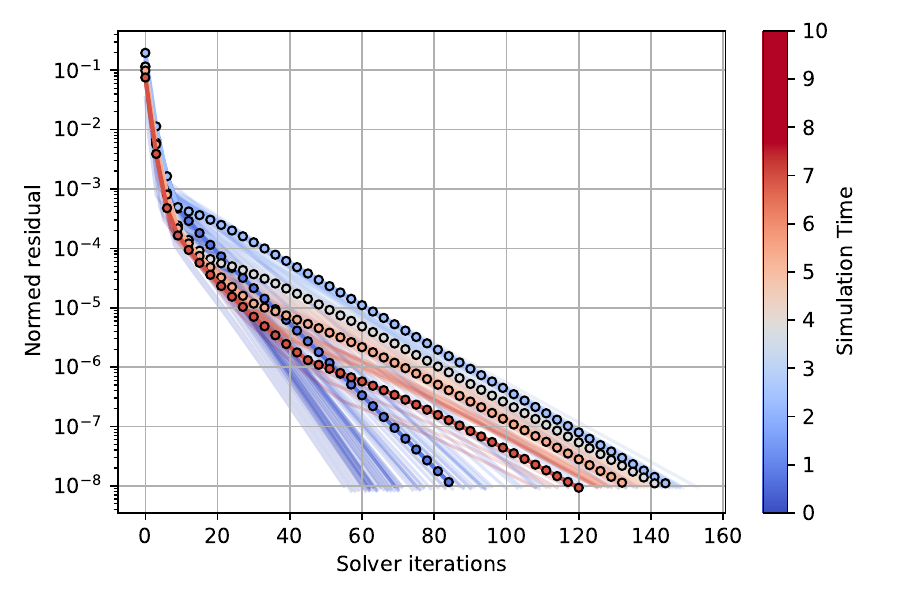}
    \caption{Mode I crack growth}
  \end{subfigure}\hfill
  \begin{subfigure}{0.5\linewidth}
    \includegraphics[height=5.5cm]{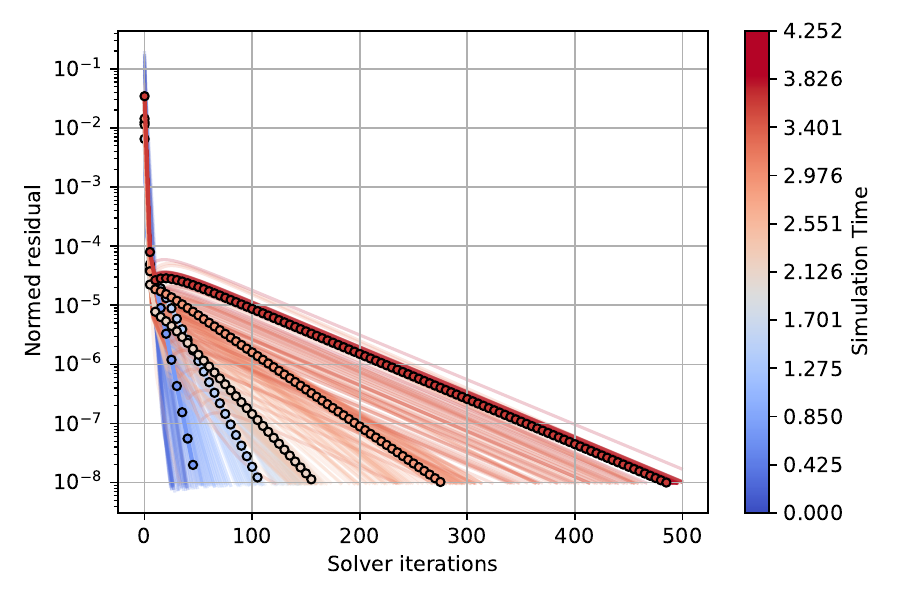}
    \caption{Stress-concentration fracture}
  \end{subfigure}
  
  \begin{subfigure}{0.5\linewidth}
    \includegraphics[height=5.5cm]{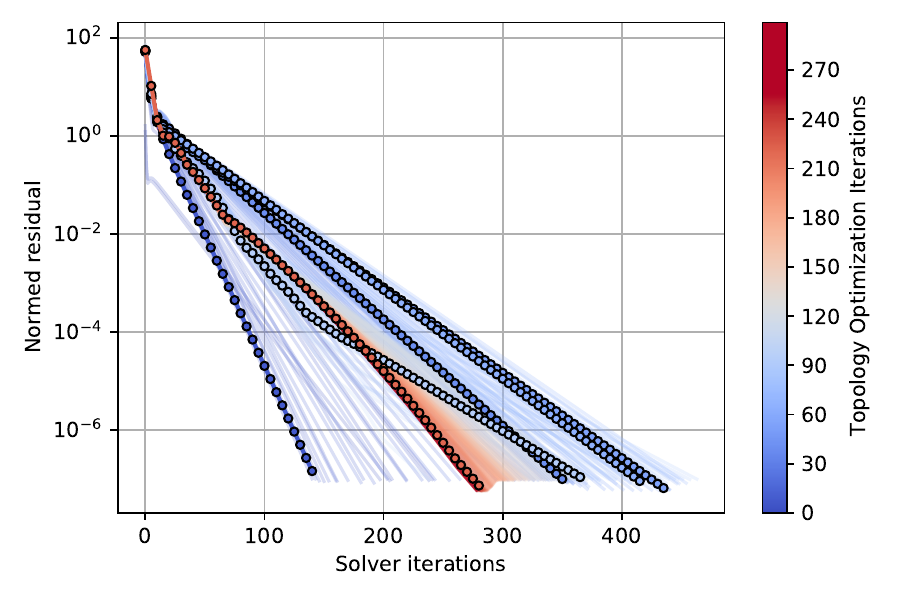}
    \caption{Topology optimization}
  \end{subfigure}
  \caption{
    Convergence of the linear elastic solver for mode I fracture, stress-concentration fracture, and topology optimization.
      In each plot, each line represents the residual vs. iteration number for a single elastic solve.
      The colors of the lines correspond to the time at which the solve occurred in the larger simulation, as well as the color in the corresponding convergence plots.
      The highlighted lines (with markers) correspond to the specific simulation points as highlighted in \cref{fig:fracture_modeI,fig:crack_lshape,fig:topop_designevolution}.
  }
  \label{fig:convergence_plots}
\end{figure}

The linear elastic, strong-form near-singular multigrid solver exhibited nearly universally linear convergence in the example problems presented in this work.
Here, we present a more detailed exposition of the solver behavior in time for the fracture cases and the topology optimization example (\cref{fig:convergence_plots}).
We note that in all cases, the results from the initial solve are not included because it is for the initial, non-regularized version of the problem.
In all examples, we observe a rapid convergence during the first 3-4 iterations.
This rapid convergence lasts until the error is reduced to $10^{-3}$, $10^{-4}$, and $10^{-1}$ for mode I, stress-concentration, and topology optimization, respectively.
A sharp reduction follows this in the convergence rate, which is almost always non-decreasing during the remaining solve.

In mode-I fracture and topology optimization, interestingly, the worst convergence occurs during the middle of the simulation; subsequently, convergence improves and approaches a constant rate.
Both exhibit a couple of solves in which the convergence rate turned sharply from a higher to a lower value; in both cases, the convergence might be described as ``piecewise linear.''

In the stress-concentration fracture case, linear convergence is constantly observed, but the convergence rate decreases steadily as the simulation progresses.
We attribute this to the increasing irregularity of the problem, resulting from large chunks of material that have been degraded, and regions of the boundary that are under-regularized. 
In other words, it appears to be the fracture model, not the solver, that is responsible for the degrading convergence.

\bibliographystyle{ieeetr} 
\bibliography{library}

\end{document}